\input ./style/arxiv-general.cfg
\documentclass[MSNbibl,number,citesort,dvips]{arxbj}
\makeatletter
   \@ifpackageloaded{graphicx}{}{\usepackage{graphicx}}
\makeatother

\volume{23}
\issue{2}
\pubyear{2017}
\firstpage{789}
\lastpage{824}
\doi{10.3150/14-BEJ679}
\docsubty{FLA}

\makeatletter
\renewcommand{\mid}{\vert}
\newcommand{\rrvert}{\vert}
\newcommand{\rrVert}{\Vert}
\newcommand{\llvert}{\vert}
\newcommand{\llVert}{\Vert}
\newtheorem{theorem}{Theorem}
\newtheorem{lemma}[theorem]{Lemma}

\newproclaim{definition}{Definition}
\newremark{rems}{Remarks}
\newremark{rem}{Remark}
\newcommand{\R}{\mathbb{R}}
\newcommand{\reals}{\mathbb{R}}
\newcommand{\En}{\mathbb{E}} 
\newcommand{\Prob}{\mathbb{P}}
\newcommand{\e}{\mathbf{e}}
\newcommand{\cC}{\mathcal{C}}
\newcommand{\cE}{\mathcal{E}}
\newcommand{\F}{\mathcal{F}}
\newcommand{\G}{\mathcal{G}}
\newcommand{\cL}{\mathcal{L}}
\newcommand{\cN}{\mathcal{N}}
\newcommand{\cP}{\mathcal{P}}
\newcommand{\X}{\mathcal{X}}
\newcommand{\Y}{\mathcal{Y}}
\newcommand{\Z}{\mathcal{Z}}
\newcommand{\cH}{\mathcal{H}}
\newcommand{\M}{\mathcal{M}}
\newcommand{\Rad}{\mathfrak{R}}
\renewcommand\P{\mathbf{P}}
\makeatother

\begin{document}
\begin{frontmatter}

\title{Empirical entropy, minimax regret and minimax risk}
\runtitle{Empirical entropy, minimax regret and minimax risk}

\begin{aug}
\author[A]{\inits{A.}\fnms{Alexander}~\snm{Rakhlin}\corref{}\thanksref{A}\ead[label=e1]{rakhlin@wharton.upenn.edu}},
\author[B]{\inits{K.}\fnms{Karthik}~\snm{Sridharan}\thanksref{B}}
\and
\author[C]{\inits{A.B.}\fnms{Alexandre B.}~\snm{Tsybakov}\thanksref{C}\ead[label=]{}}
\address[A]{Department of Statistics, University of Pennsylvania, Philadelphia, PA 19104, USA.\\ \printead{e1}}
\address[B]{Department of Computer Science, Cornell University, Ithaca, NY 14853, USA}
\address[C]{Laboratoire de Statistique, CREST-ENSAE, 92245 Malakoff Cedex, France}
\end{aug}

%
\received{\smonth{3} \syear{2014}}
%
\revised{\smonth{7} \syear{2014}}

%
\begin{abstract}
We consider the random design regression model with square loss. We
propose a method that aggregates empirical minimizers (ERM) over
appropriately chosen random subsets and reduces to ERM in the extreme
case, and we establish sharp oracle inequalities for its risk. We show
that, under the $\varepsilon^{-p}$ growth of the empirical $\varepsilon
$-entropy,\vspace*{1pt} the excess risk of the proposed method attains the rate
$n^{-2/(2+p)}$ for $p\in(0,2)$ and $n^{-1/p}$ for $p> 2$ where
$n$ is the sample size.
Furthermore, for $p\in(0,2)$, the excess risk rate matches the
behavior of the minimax risk of function estimation in regression
problems under the well-specified model. This yields a conclusion that
the rates of statistical estimation in well-specified models (minimax
risk) and in misspecified models (minimax regret) are equivalent in the
regime $p\in(0,2)$. In other words, for $p\in(0,2)$ the problem of
statistical learning enjoys the same minimax rate as the problem of
statistical estimation. On the contrary, for $p>2$ we show that the
rates of the minimax regret are, in general, slower than for the
minimax risk. Our oracle inequalities also imply the $v\log(n/v)/n$
rates for Vapnik--Chervonenkis type classes of dimension $v$ without the
usual convexity assumption on the class; we show that these rates are
optimal. Finally, for a slightly modified method, we derive a bound on
the excess risk of $s$-sparse convex aggregation improving that of
Lounici [\textit{Math. Methods Statist.} \textbf{16} (2007) 246--259] and providing the optimal rate.
\end{abstract}

%
\end{frontmatter}

\section{Introduction}

Let $D_n=\{(X_1,Y_1),\ldots,(X_n,Y_n)\}$ be an i.i.d. sample from
distribution $P_{XY}$ of a pair of random variables $(X,Y)$, $X\in\X
$, $Y\in\Y$ where $\X$ is any set and $\Y$ is a subset of $\R$. We
consider the problem of prediction of $Y$ given $X$. For any measurable
function $f\dvtx\X\to\Y$ called the predictor, we define the prediction
risk under squared loss:
\[
L(f)=\En_{XY} \bigl[\bigl(f(X)-Y\bigr)^2\bigr], %
\]
where $\En_{XY}$ is the expectation with respect to $P_{XY}$.
Let now $\F$ be a class of functions from $\X$ to~$\Y$ and assume
that the aim is to mimic the best predictor in this class.
This means that we want to find an estimator $\hat{f}$ based on the
sample $D_n$ and having a small excess risk
%
\begin{eqnarray}
\label{eqexcriskdef} L(\hat{f})-\inf_{f\in\F} L(f)
\end{eqnarray}
in expectation or with high probability. The minimizer of $L(f)$ over
all measurable functions is the regression function $\eta(x) = \En
_{XY}[Y\mid X=x]$ and it is straightforward to see that for the expected
excess risk we have
%
\begin{eqnarray}
\label{0} \cE_{\F}(\hat{f}) \triangleq\En L(\hat{f})-\inf
_{f\in\F} L(f) = \En\llVert\hat{f}-\eta\rrVert^2 -
\inf_{f\in\F} \llVert f-\eta\rrVert^2,
\end{eqnarray}
where $\En$ is the generic expectation sign, $\llVert f\rrVert
^2={\int f^2(x)
P_X(\mathrm{d}x)}$, and $P_X$ denotes the marginal distribution of $X$. The
left-hand side of (\ref{0}) has been studied within Statistical Learning Theory characterizing the error of ``agnostic learning'' \cite
{vapnik1974theory,DevrGyorLug96,koltchinskii2011oracle}, while the object on the right-hand side has
been the topic of oracle inequalities in nonparametric statistics \cite
{nemirovski2000topics,Tsy09}, and in the literature on
aggregation \cite{tsybakov03optimal,rigollet2011exponential}.
Upper bounds on the right-hand side of (\ref{0}) are called \emph
{sharp} oracle inequalities, which refers to constant $1$ in front of
the infimum over $\F$. However, some of the key results in the
literature were only obtained with a constant greater than $1$, that
is, they yield upper bounds for the difference
%
\begin{eqnarray}
\label{eqnonexactoracle} \En\llVert\hat{f}-\eta\rrVert^2 - C \inf
_{f\in\F} \llVert f-\eta\rrVert^2
\end{eqnarray}
with $C>1$ and not for the excess risk. In this paper, we obtain sharp
oracle inequalities, which allows us to consider the excess risk
formulation of the problem as described above.

In what follows we assume that $\Y=[0,1]$. For results in expectation,
the extension to unbounded $\Y$ with some condition on the tails of
the distribution is straightforward. For high probability statements,
more care has to be taken, and the requirements on the tail behavior
are more stringent. To avoid this extra level of complication, we
assume boundedness.

From a minimax point of view, the object of interest in Statistical
Learning Theory can be written as the \textit{minimax regret}
%
\begin{eqnarray}
\label{eqminimaxregret} V_n(\F) = \inf_{\hat{f}}\sup
_{P_{XY}\in\cP} \Bigl\{ \En L(\hat{f})-\inf_{f\in\F} L(f)
\Bigr\},
\end{eqnarray}
where $\cP$ is the set of all probability distributions on $\X\times
\Y$ and $\inf_{\hat{f}}$ denotes the infimum over all estimators. We
observe that the study of this object leads to a \emph
{distribution-free} theory, as no model is assumed. Instead, the goal
is to achieve predictive performance competitive with a reference class
$\F$. In view of (\ref{0}), an equivalent way to write $V_n(\F)$ is
%
\begin{eqnarray}
\label{eqminimaxregretnorms} V_n(\F) = \inf_{\hat{f}}\sup
_{P_{XY}\in\cP} \Bigl\{ \En\llVert\hat{f}-\eta\rrVert^2 -
\inf_{f\in\F} \llVert f-\eta\rrVert^2 \Bigr\}.
\end{eqnarray}
The minimax regret can be interpreted as a measure of performance of
estimators for misspecified models. The study of $V_n(\F)$ will be
further referred to as \textit{misspecified model} setting.

A special instance of the minimax regret has been studied in the
context aggregation of estimators, with the aim to characterize optimal
rates of aggregation, cf., for example,
\cite{tsybakov03optimal,rigollet2011exponential}. There, $\F$~is a
subclass of the linear span of $M$ given functions $f_1,\ldots,f_M$,
for example, their convex hull or sparse linear (convex) hull.
Functions $f_1,\ldots,f_M$ are interpreted as some initial estimators
of the regression function $\eta$ based on another sample from the
distribution of $(X,Y)$. This sample is supposed to be independent from
$D_n$ and is considered as frozen when dealing with the minimax regret.
The aim of aggregation is to construct an estimator $\hat f$, called
the aggregate, that mimics the best linear combination of $f_1,\ldots
,f_M$ with coefficients of the combination lying in a given set in $\R
^M$. Our results below apply to this setting as well. We will provide
their consequences for some important examples of aggregation.

In the standard nonparametric regression setting, it is assumed that
the model is \textit{well-specified}, that is, we have $Y_i=f(X_i)+\xi_i$
where the random errors $\xi_i$ satisfy $\En(\xi_i\mid X_i)=0$ and $f$
belongs to a given functional class $\F$. Then $f=\eta$ and the
infimum on the right-hand side of (\ref{0}) is zero. The value of
reference characterizing the best estimation in this problem is the
\textit{minimax risk}
%
\begin{eqnarray}
\label{eqminimaxestimate} W_n(\F) = \inf_{\hat{f}}\sup
_{P_{XY}\in\mathcal{P}_\F} \En\llVert\hat{f}-\eta\rrVert^2,
\end{eqnarray}
where $\mathcal{P}_\F$ is the set of all distributions $P_{XY}$ on
$\X\times\Y$ such that $\eta\in\F$.
It is not difficult to see that
\[
W_n(\F) \leq V_n(\F) , %
\]
yet the minimax risk and the minimax regret are quite different and the
question is whether the two quantities can be of the same order of
magnitude for particular $\F$. We show below that the answer is
positive 
for major cases of interest except for very massive classes $\F$,
namely, those having the empirical $\varepsilon$-entropy of the order
$\varepsilon^{-p}$, $p>2$, for small $\varepsilon$.
We also prove that this entropy condition is tight in the sense that
the minimax regret and the minimax risk can have different rates of
convergence when it is violated. Furthermore, we show that the optimal
rates for the minimax regret and minimax risk are attained by one and
the same procedure -- the aggregation-of-leaders estimator -- that we
introduce below.

Observe a certain duality between $W_n(\F)$ and $V_n(\F)$. In the
former, the assumption about the reality is placed on the way data are
generated. In the latter, no such assumption is made, yet the
assumption is placed in the term that is being subtracted off. As we
describe in Section~\ref{sechistorical}, the study of these two
quantities represents two parallel developments: the former has been a
subject mostly studied within nonparametric statistics, while the
second -- within Statistical Learning Theory. We aim to bring out a
connection between these two objects. In Section~\ref{secapprox}, we
introduce a more general risk measure that realizes a smooth transition
between $W_n(\F)$ and $V_n(\F)$ depending on the magnitude of the
approximation error. The minimax risk and the minimax regret appear as
the two extremes of this scale.

The paper is organized as follows. In Section~\ref{secmain}, we
present the aggregation-of-leaders estimator and the upper bounds on
its risk. These include the main oracle inequality in Theorem~\ref
{lemmain} and its consequences for particular classes $\F$ in
Theorems~\ref{thmmainc}--\ref{thmmainaggr}.
Section~\ref{secapprox} discusses a more general setting allowing for
a smooth transition between $W_n(\F)$ and $V_n(\F)$ in terms of the
approximation error. Lower bounds for the minimax risk and minimax
regret are proved in Section~\ref{seclower}.
In Section~\ref{seccompare}, we compare the aggregation-of-leaders
estimator with the two closest competitors -- skeleton aggregation and
global ERM.
Section~\ref{sechistorical} provides an overview and comparison of
our results to those in the literature. Proofs of the theorems are
given in Sections~\ref{secproofsthms}--\ref{secproofslower}. The
\hyperref[secappendix]{Appendix} contains some technical results and proofs of the lemmas.

\section{Notation}\label{secnotation}

Set $\Z=\X\times\Y$. For $S=\{z_1,\ldots,z_n\}\in\Z^n$ and a
class $\G$ of real-valued functions on $\Z$, consider the Rademacher
average of $\G$:
\[
\hat{\Rad}_n(\G, S) = \En_{\sigma} \Biggl[ \sup
_{g\in\G}\frac
{1}{n}\sum_{i=1}^n
\sigma_i g(z_i) \Biggr],
\]
where $\En_{\sigma}$ denotes the expectation with respect to the
joint distribution of i.i.d. random variables $\sigma_1,\ldots,\sigma_n$
taking values $1$ and $-$1 with probabilities $1/2$. Let
\[
\Rad_n(\G) = \sup_{S\in\Z^n} \hat{\Rad}_n(
\G,S).
\]
Given $r>0$, we denote by $\G[r,S]$ the set of functions in $\G$ with
empirical average at most $r$ on~$S$:
\[
\G[r,S] = \Biggl\{g\in\G\dvt \frac{1}{n}\sum_{i=1}^n
g(z_i) \leq r \Biggr\}.
\]

Any function $\phi_n\dvtx[0,\infty)\mapsto\reals$ satisfying
%
\begin{eqnarray}
\label{eqphindef} \sup_{S\in\Z^n}\hat{\Rad}_n \bigl(
\G[r,S], S \bigr) \leq\phi_n(r)
\end{eqnarray}
for all $r>0$ will be called an upper function for the class $\G$. We
will sometimes write $\phi_n(r)=\phi_n(r,\G)$ to emphasize the
dependence on $\G$. It can be shown (cf., e.g., Lemma~\ref
{lemboundonr-star} below) that any class of uniformly bounded
functions admits an upper function satisfying the sub-root property:
$\phi_n$ is non-negative, non-decreasing, and $\phi_n(r)/\sqrt{r}$
is non-increasing. We will denote by $r^*=r^*(\G)$ the corresponding
\textit{localization radius}, that is, an upper bound on the largest
solution of the equation $\phi_n(r)=r$. Clearly, $r^*$ is not uniquely
defined since we deal here with upper bounds.

We write $\ell\circ f$ for the function $(x,y)\mapsto(f(x)-y)^2$ and
$\ell\circ\F$ for the class of functions $ \{\ell\circ f\dvt f\in
\F\}$. Thus,
\[
(\ell\circ\F)[r,S] = \Biggl\{ \ell\circ f\dvt f\in\F, \frac
{1}{n}\sum
_{i=1}^n (\ell\circ f)
(x_i,y_i) \leq r \Biggr\}
\]
for $S=\{z_1,\ldots,z_n\}$ with $z_i=(x_i,y_i)$.

For any bounded measurable function $g\dvtx\Z\to\R$, we set $Pg= \En
g(Z)$, where $Z=(X,Y)$, and $P_ng=\frac{1}{n} \sum_{i=1}^n g(Z_i)$,
where $Z_i=(X_i,Y_i)$.
For $S=\{z_1,\ldots,z_n\}\in\Z^n$ with $z_i=(x_i,y_i)$ consider the
empirical $\ell_2$ pseudo-metric
\[
d_S (f,g)= \Biggl(\frac{1}{n} \sum
_{i=1}^n \bigl\llvert f(x_i)-g(x_i)
\bigr\rrvert^2 \Biggr)^{1/2} , %
\]
and for any $\varepsilon>0$ denote by $\cN_2(\F, \varepsilon, S)$ the
{$\varepsilon$-covering number} of a class $\F$ of real-valued functions
on $\X$ with respect to this pseudo-metric. Recall that a covering
number at scale $\varepsilon$ is the smallest number of balls of radius
$\varepsilon$ required to cover the set. Denote by $\cN_\infty(\F,
\varepsilon, S)$ the $\varepsilon$-covering number of the class $\F$ with
respect to the supremum norm (over $S$).

Although not discussed here explicitly, some standard measurability
conditions are needed to apply results from the theory of empirical
processes as well as to ensure that the ERM estimators we consider
below are measurable. This can be done in a very general framework and
we assume throughout that these conditions are satisfied. For more
details we refer to Chapter~5 of \cite{Dudley99}, see also \cite
{koltchinskii2011oracle}, page 17.

The minimum risk on the class of functions $\F$ is denoted by
\[
L^*=\inf_{f\in\F} L(f).
\]
%
Let $\lceil x \rceil$ denote the minimal integer strictly greater than
$x\in\R$, and $\llvert\F\rrvert$ the cardinality of $\F
$. Notation $C$ will be
used for positive constants that can vary on different occasions; these
are absolute constants unless their dependence on some parameters is
explicitly mentioned. We will also assume throughout that $n\ge5$.

\section{Main results}\label{secmain}

In this section, we introduce the estimator studied along the paper,
state the main oracle inequality for its risk and provide
corollaries for the minimax risk and minimax regret. The estimation
procedure comprises three steps. The first step is to construct a
random $\varepsilon$-net on $\F$ with respect to the empirical $\ell_2$
pseudo-metric and to form the induced partition of $\F$. The second
step is to compute empirical risk minimizers (in our case, the least
squares estimators) over cells of this random partition. Finally,
the third step is to aggregate these minimizers using a suitable
aggregation procedure. If the radius $\varepsilon$ of the initial net
is taken to be large enough, the method reduces to the global
empirical risk minimization (ERM) over the class $\F$. While the
global ERM is, in general, suboptimal (cf. the discussion in
Sections~\ref{seccompare} and \ref{sechistorical} below), the
proposed method enjoys the optimal rates. We call our method the
aggregation-of-leaders procedure since it aggregates the best
solutions obtained in cells of the partition.

To ease the notation, assume that we have a sample $D_{3n}$ of size
$3n$ and we divide it into three parts: $D_{3n}=S\cup S' \cup S''$,
where the subsamples $S, S' , S''$ are each of size $n$. Fix
$\varepsilon
>0$. Let $d_S(f,g)$
be the empirical $\ell_2$ pseudo-metric associated with the subsample
$S$ of cardinality $n$, and
\[
N = \cN_2(\F, \varepsilon, S).
\]
Clearly, $N$ is finite since $\F$ is included in the set of all
functions with values in $[0,1]$, which is totally bounded with respect
to $d_S(\cdot,\cdot)$.
Let $\hat{c}_1,\ldots,\hat{c}_N$ be an $\varepsilon$-net on $\F$ with
respect to $d_S(\cdot,\cdot)$. We assume without loss of generality
that it is \emph{proper}, that is, $\hat{c}_i\in\F$ for $i=1,\ldots
,N$, and that $N\ge2$. Let $\hat{\F}_1^S,\ldots,\hat{\F}_N^S$ be
the following partition of $\F$ induced by $\hat{c}_i$'s:
\[
\hat{\F}_i^S = \hat{\F}_i^S (
\varepsilon)= \Bigl\{f\in\F\dvt i \in\mathop{\operatorname{argmin}}_{j=1,\ldots,N}
d_S(f, \hat{c}_j) \Bigr\}
\]
with ties broken in an arbitrary way. Now, for each $\hat{\F}_i^S$,
define the least squares estimators over the subsets $\hat{\F}_i^S$
with respect to the second subsample $S'$:
%
\begin{equation}
\label{defest} \hat{f}_i^{S,S'} \in\mathop{\operatorname{argmin}}_{f\in\hat{\F}_i^S} \frac
{1}{n}\sum_{(x,y)\in S'}
\bigl(f(x)-y\bigr)^2.
\end{equation}
We will assume that such a minimizer exists; a simple modification of
the results is
possible if $\hat{f}_i^{S,S'}$ is an approximate solution of (\ref{defest}).

Finally, at the third step we use the subsample $S''$ to aggregate the
estimators $ \{\hat{f}_1^{S,S'}, \ldots,\break \hat{f}_N^{S,S'}
\}$.
We call a function $\tilde f (x, D_{3n})$ with values in $\Y$ a \textit
{sharp $\mathrm{MS}$-aggregate}\footnote{Here, $\mathrm{MS}$-aggregate is an abbreviation
for \textit{model selection type aggregate}. The word \textit{sharp}
indicates that (\ref{eqaggregationbound}) is an oracle inequality
with leading constant 1.} if it has the following property.

\begin{SMS*}
There exists a constant $C>0$
such that, for any $\delta>0$,
%
\begin{eqnarray}
\label{eqaggregationbound} L (\tilde{f} ) \leq\min_{i=1,\ldots,N} L
\bigl(\hat
{f}_i^{S,S'} \bigr) + C\frac{\log(N/\delta) }{n}
\end{eqnarray}
with probability at least $1-\delta$ over the sample $S''$,
conditionally on $S\cup S'$.
\end{SMS*}

Note that, in (\ref{eqaggregationbound}), the subsamples $S, S'$ are
fixed, so that the estimators $\hat{f}_i^{S,S'}\triangleq g_i$ can be
considered as fixed (non-random) functions, and $\tilde f$ as a
function of $S''$ only. There exist several examples of sharp
$\mathrm{MS}$-aggregates of fixed functions $g_1,\ldots,g_N$ \cite
{audibert2007progressive}, page~5,
\cite{lecue2009aggregation}, Theorem~\textup{B}, \cite{LecueRigollet},
Theorem~\textup{A}. They are realized as mixtures:
%
\begin{equation}
\label{mixture} \tilde{f} = \sum_{i=1}^N
\theta_i g_i = \sum_{i=1}^N
\theta_i \hat{f}_i^{S,S'},
\end{equation}
where $\theta_i$ are some random weights measurable with respect to
$S''$. Either of the aggregates of \cite
{audibert2007progressive,lecue2009aggregation,LecueRigollet} satisfy
the sharp MS-aggregation property and thus can be used at the third
step of our procedure.

%
\begin{definition}
We call an \emph{aggregation-of-leaders} estimator any estimator
$\tilde{f}$ defined by the above three-stage procedure with sharp
$\mathrm{MS}$-aggregation at the third step.
\end{definition}

The next theorem provides the main oracle inequality for
aggregation-of-leaders estimators.

%
\begin{theorem}
\label{lemmain}
Let $\Y=[0,1]$ and $0\leq f\leq1$ for all $f\in\F$. Let $r^*=r^*(\G
)$ denote a localization radius of $\G=\{(f-g)^2\dvt f,g\in\F\}$.
Consider an aggregation-of-leaders estimator $\tilde{f}$ defined by
the above three-stage procedure. Then there exists an absolute constant
$C>0$ such that for any $\delta>0$, with probability at least
$1-2\delta$,
%
\begin{eqnarray}
\label{1} L(\tilde{f}) \le\inf_{f\in\F} L(f) + C \biggl(
\frac{\log(\cN
_2(\F,\varepsilon,S)/\delta)}{n}+ \Xi\bigl(n, \varepsilon, S'\bigr)
\biggr),
\end{eqnarray}
where
%
\begin{eqnarray}
\label{2} \Xi\bigl(n, \varepsilon, S'\bigr)&=&\gamma
\sqrt{r^*} + \inf_{\alpha\geq0} \biggl\{ \alpha+ \frac
{1}{\sqrt{n}}\int
_{\alpha}^{C\gamma} \sqrt{\log\cN_2\bigl(\F,
\rho, S'\bigr)} \,\mathrm{d}\rho\biggr\}
\end{eqnarray}
with $\gamma= \sqrt{\varepsilon^2 + r^* + \beta}$ and $\beta= (\log
(1/\delta) + \log\log n)/n$.

\end{theorem}

\begin{rems*}
\begin{enumerate}
\item[1.] The term $\Xi(n, \varepsilon, S')$ in Theorem~\ref{lemmain} is a
bound on the rate of convergence of the excess risk of ERM $\hat
{f}_{i}^{S,S'}$ over the cell $\hat{\F}^{S}_{i}$.
If, in particular instances, there exists a sharper bound for the rate
of ERM, one can readily use this bound instead of the expression for
$\Xi(n, \varepsilon, S')$ given in Theorem~\ref{lemmain}. 
%
\item[2.] The partition with cells $\hat{\F}^{S}_{i}$ defined above can
be viewed as a default option. In some situations, we may better tailor
the (possibly overcomplete) partition to the geometry of~$\F$. For
instance, in the aggregation context (cf. Theorem~\ref{thmmainaggr}
below), $\F$ is union of convex sets. We choose each convex set as an
element of the partition, and use the rate for ERM over individual
convex sets instead of the overall rate $\Xi(n, \varepsilon, S')$. In
this case, the partition is non-random. Another example, when $\F$ is
isomorphic to a subset of $\R^M$, is a partition of~$\R^M$ into a
union of linear subspaces of all possible dimensions. In this case, the
``cells'' are linear subspaces and aggregating the least squares
estimators over cells is analogous to sparsity pattern aggregation
considered in \cite{rigollet2011exponential,RigTsySTS12}.

\item[3.] In Theorem~\ref{lemmain} we can use the localization radius
$r^*=r^*(\hat{\G}_i)$ for $\hat{\G}_i=\{(f-g)^2\dvt f,g\in\hat{\F
}^S_i\}$ instead of the larger quantity $r^*(\G)$. 
Inspection of the proof shows that the oracle inequality~(\ref{1})
generalizes to
%
\begin{eqnarray}
\label{1a} L(\tilde{f}) \le\min_{i=1,\ldots,N}\inf
_{f\in\hat{\F}^{S}_i} \bigl\{L(f) + C \bigl(\beta+ \Xi_i\bigl(n,
\varepsilon, S'\bigr)\bigr) \bigr\},
\end{eqnarray}
where $\Xi_i(n, \varepsilon, S')$ is defined in the same way as $\Xi(n,
\varepsilon, S')$ with the only difference that $r^*(\G)$ is replaced by
$r^*(\hat{\G}_i)$.
\end{enumerate}
\end{rems*}

The oracle inequality (\ref{1}) of Theorem~\ref{lemmain} depends on
two quantities that should be specified: the entropy $ \log\cN_2(\F,
\cdot, \cdot)$, and the localization radius $r^*$. The crucial role
in determining the rate belongs to the empirical entropies. We further
replace in~(\ref{1}) these random entropies by their upper bound
\[
\cH_2(\F,\rho)= \sup_{S\in\Z^n}\log\cN_2(
\F, \rho, {S}), %
\]
and refer to the above quantity as the \emph{empirical entropy}.

The next theorem is a corollary of Theorem~\ref{lemmain} in the
case of polynomial growth of the empirical entropy characteristic
for nonparametric estimation problems. It gives upper bounds on the
minimax regret and on the minimax risk.
%
\begin{theorem}
\label{thmmainc}
Let $\Y=[0,1]$ and $\cH_2(\F,\rho)\le A\rho^{-p}$, $\forall\rho
>0$, for some constants $A<\infty$, $p>0$. Let $\tilde{f}$ be an
aggregation-of-leaders estimator defined by the above three-stage
procedure with the covering radius $\varepsilon= n^{-1/(2+p)}$.
There exist constants $C_p>0$ depending only on $A$ and $p$ such that:
\begin{longlist}[(ii)]
\item[(i)] Let $0\leq f\leq1$ for all $f\in\F$. For the estimator
$\tilde{f}$ we have: 
%
\begin{eqnarray}
\label{4} V_n(\F)&\le&\sup_{P_{XY}\in\cP} \Bigl\{ \En
\llVert\tilde{f}-\eta\rrVert^2 -\inf_{f\in\F}
\llVert f-\eta\rrVert^2 \Bigr\}
\nonumber\\[-8pt]\\[-8pt]\nonumber
&\le& \cases{ \displaystyle
C_p n^{-2/(2+p)}, &\quad if $p\in(0,2)$,
\vspace*{3pt}\cr
C_p
n^{-1/2}\log(n), &\quad if $p=2$,
\vspace*{3pt}\cr
C_p n^{-1/p}, &
\quad if $p\in(2,\infty)$.}
\end{eqnarray}
\item[(ii)] When the model is well-specified,
then for the estimator $\tilde{f}$ we have:
%
\begin{eqnarray}
\label{5} W_n(\F)\le\sup_{P_{XY}\in{\mathcal P}_\F} \En\llVert
\tilde{f}-\eta\rrVert^2 \le C_p n^{-2/(2+p)} \qquad
\forall p>0.
\end{eqnarray}
\end{longlist}
\end{theorem}

The proof of Theorem~\ref{thmmainc} is given in
Section~\ref{secproofsthms}. The first conclusion of this
theorem is that the minimax risk $W_n(\F)$ has the same rate of
convergence as the minimax regret $V_n(\F)$ for $p\in(0,2)$. For
example, if $\F$ is a class of functions on $\R^d$ with bounded
derivative of order $k$, the entropy bound required in the
theorem holds with exponent $p=d/k$, as follows from \cite
{kolmogorov1959}. In
this case, Theorem~\ref{thmmainc} yields that, for $k\ge d/2$,
both $W_n(\F)$ and $V_n(\F)$ converge with the usual nonparametric
rate $n^{-{2k}/(2k +d)}$ while for $k<d/2$
(corresponding to very irregular functions) the rate of the minimax
regret deteriorates to $n^{-k/d}$. In Section~\ref{seclower}, we will
show that the bounds of Theorem~\ref{thmmainc} for $p< 2$
are tight in the sense that there exists a marginal distribution of $X$
and a class $\F$ of regression functions satisfying the above entropy
assumptions such that the bounds (\ref{4}) and (\ref{5}) cannot be
improved for $p< 2$.

The second message of Theorem~\ref{thmmainc} is that $W_n(\F)$ has
faster rate than $V_n(\F)$ for $p>2$, that is, for very massive
classes $\F$. Note that here we compare only the upper bounds.
However, in Section~\ref{seclower} we will provide a lower bound
showing that the effect indeed occurs. Namely, we will exhibit a
marginal distribution of $X$ and a class $\F$ of regression functions
satisfying the above entropy assumptions such that $V_n(\F)$ is of the
order $n^{-1/(p-1)}$, which is slower than the rate $n^{-2/(2+p)}$ for $W_n(\F)$.

Observe also that in both cases, $p\in(0,2)$ and $p\in[2, \infty)$,
we can use the same value $\varepsilon= n^{-1/(2+p)}$ to obtain the rates
given in (\ref{4}). We remark that this $\varepsilon$ satisfies the
balance relation
\[
n\varepsilon^2 \asymp\cH_2(\F,\varepsilon). %
\]
We will further comment on this choice in Section~\ref{seccompare}.

We now turn to the consequences of Theorem~\ref{lemmain} for low
complexity classes $\F$, such as Vapnik--Chervonenkis (VC) classes and
intersections of balls in finite-dimensional spaces. They roughly
correspond to the case ``$p\approx0$'', and the rates for the minimax
risk $W_n(\F)$ are the same as for the minimax regret $V_n(\F)$.

Assume first that the empirical covering numbers of $\F$ exhibit the growth
%
\begin{equation}
\label{eq1thmmainVC} \sup_{S\in\Z^n}\cN_2(\F,\rho,S)\le(A/
\rho)^{v},
\end{equation}
$\forall\rho>0$, with some constants $A<\infty$, $v>0$.
Such classes $\F$ are called VC-\textit{type classes} with
VC-dimension $\operatorname{VC}(\F)=v$. We will also call them \textit{parametric
classes}, as opposed to nonparametric classes considered in
Theorem~\ref{thmmainc}. Indeed, entropy bounds as in (\ref
{eq1thmmainVC}) are associated to compact subsets of $v$-dimensional
Euclidean space. Other example is given by the VC-subgraph classes with
VC-dimension $v$, that is, classes of functions $f$ whose subgraphs
$C_f = \{(x,t)\in\X\times\R\dvt f(x)\geq t\}$ form a
Vapnik--Chervonenkis class with VC-dimension $v$.

%
\begin{theorem}[(Bounds for VC-type classes)]\label{thmmainVC}
Assume that $\Y=[0,1]$ and the empirical covering numbers satisfy
(\ref{eq1thmmainVC}).
Let $0\leq f\leq1$ for all $f\in\F$, and let $\tilde{f}$ be an
aggregation-of-leaders estimator defined by the above three-stage
procedure with $\varepsilon= n^{-1/2}$. If $n\ge C_Av$ for a
large enough constant $C_A>1$ depending only on $A$, then there exists
a constant $C>0$ depending only on $A$ such that
%
\begin{eqnarray}
\label{eq2thmmainVC} V_n(\F)\le\sup_{P_{XY}\in\cP} \Bigl\{ \En
\llVert\tilde{f}-\eta\rrVert^2 -\inf_{f\in\F}
\llVert f-\eta\rrVert^2 \Bigr\}\le C \frac{v}{n} \log\biggl(
\frac{en}{v} \biggr) .
\end{eqnarray}
\end{theorem}

The rate of convergence of the excess risk as in (\ref
{eq2thmmainVC}) for VC-type classes has been obtained previously
under the assumption that $L^*=0$ or for convex classes $\F$ (see
discussion in Section~\ref{sechistorical} below). Theorem~\ref
{thmmainVC} does not rely on either of these assumptions.

In Section~\ref{seclower}, we show that the bound of Theorem~\ref
{thmmainVC} is tight; there exists a function class such that, for
any estimator, there exists a distribution on which the estimator
differs from the regression function by at least $C(v/n)\log(en/v)$
with positive fixed probability. So, the extra logarithmic factor $\log
(en/v)$ in the rate is necessary, even when the model is well-specified.

The next theorem deals with classes of functions
\[
\F= \F_\Theta\triangleq\Biggl\{\mathsf{f}_\theta= \sum
_{i=1}^M \theta_j f_j \dvt
\theta=(\theta_1,\ldots,\theta_M) \in\Theta\Biggr\},
\]
where $\{f_1,\ldots,f_M\}$ is a given collection of $M$ functions on
$\X$ with values in $\Y$, and $\Theta\subseteq\R^M$ is a given set
of possible mixing coefficients
$\theta$. Such classes arise in the context of aggregation, cf., for
example, \cite{tsybakov03optimal,rigollet2011exponential},
where the main problem is to study the behavior of the minimax regret
$V_n(\F_\Theta)$ based on the geometry of $\Theta$. For the case of
fixed rather than random design, we refer to \cite
{rigollet2011exponential} for a comprehensive treatment. Here, we deal
with the random design case and consider the sets $\Theta$ defined as
intersections of $\ell_0$-balls with the simplex. For an integer $1\le
s\le M$, the $\ell_0$-ball with radius $s$ is defined by
\[
B_0(s)=\bigl\{\theta\in\R^M\dvt \llvert\theta\rrvert
_0\le s\bigr\},
\]
where $\llvert\theta\rrvert_0$ denotes the number of
non-zero components of
$\theta$. 
We will also consider the simplex
\[
\Lambda_M= \Biggl\{\theta\in\R^M\dvt \sum
_{j=1}^M\theta_j=1,
\theta_j\ge0, j=1,\ldots,M \Biggr\}. %
\]
Then, model selection type aggregation (or $\mathrm{MS}$-\textit{aggregation})
consists in constructing an estimator $\tilde{f}$ that mimics the best
function among $f_1,\ldots, f_M$, that is, the function that attains
the minimum $\min_{j=1,\ldots,M} \llVert f_j-\eta\rrVert
^2$. In this case,
$\F_\Theta=\{f_1,\ldots, f_M\}$ or equivalently $\Theta=\Theta
^{\mathrm{MS}}\triangleq\{\e_1,\ldots,\e_M\}=\Lambda_M\cap B_0(1)
$, where $\e_1,\ldots,\e_M$ are the canonical basis vectors in $\R
^M$. \textit{Convex aggregation} (or $C$-aggregation) consists in
constructing an estimator $\tilde{f}$ that mimics the best function in
the convex hull $\F=\operatorname{conv}(f_1,\ldots, f_M)$, that is, the
function that attains the minimum $\min_{\theta\in\Lambda_M} \llVert
\mathsf{f}_\theta-\eta\rrVert^2$. In this case, $\F=\F
_\Theta$ with
$\Theta=\Theta^{\mathrm{C}}\triangleq\Lambda_M$. Finally, given an
integer $1\le s\le M$, the $s$-\textit{convex aggregation} consists in
mimicking the best convex combination of at most $s$ among the
functions $f_1,\ldots, f_M$. This corresponds to the set $\Theta
^{\mathrm{C}}(s)= \Lambda_M\cap B_0(s)$. Note that $\mathrm{MS}$-aggregation
and convex aggregation are particular cases of $s$-convex aggregation:
$\Theta^{\mathrm{MS}}=\Theta^{\mathrm{C}}(1)$ and $\Theta^{\mathrm
{C}}=\Theta^{\mathrm{C}}(M)$.

For the aggregation setting, we modify the definition of cells $\hat\F
^S_i$ as discussed in Remark~\hyperref[r2]{2}. Consider the partition $\Theta
^{\mathrm{C}}(s)=\bigcup_{m=1}^s\bigcup_{\nu\in I_m} \F_{\nu,m}$
where $I_m$ is the set of all subsets $\nu$ of $\{1,\ldots,M\}$ of
cardinality $\llvert\nu\rrvert=m$, and $\F_{\nu,m}$ is
the convex hull of
$f_j$'s with indices $j\in\nu$. We use the deterministic cells
\[
\{\F_1,\ldots,\F_N\}=\{\F_{\nu,m}, m=1,\ldots,s,
\nu\in I_m\}
\]
instead of random ones $\hat\F^S_i$. Note that the subsample $S$ is
not involved in this construction. We keep all the other ingredients of
the estimation procedure as described at the beginning of this section,
and we denote the resulting estimator $\tilde{f}$. Then, using the
subsample $S$, we complete the construction by aggregating (via a sharp
$\mathrm{MS}$-aggregation procedure) only two estimators, $\tilde{f}$ and the
least squares estimator on $\Lambda_M$. The resulting aggregate is
denoted by $\tilde{f}^*$.

%
\begin{theorem}[(Bounds for \textit{s}-convex aggregation)]\label{thmmainaggr}
Let $\Y=[0,1]$, and $0\leq f_j\leq1$ for $j=1,\ldots,M$. 
Then there exists an absolute constant $C>0$ such that
%
\begin{eqnarray}
\label{eq1thmmainaggr} V_n(\F_{\Theta^{\mathrm{C}}(s)})\le\sup
_{P_{XY}\in\cP}
\Bigl\{ \En\bigl\llVert\tilde{f}^*-\eta\bigr\rrVert^2 -\inf
_{\theta\in
\Theta^{\mathrm
{C}}(s)} \llVert\mathsf{f}_\theta-\eta\rrVert
^2 \Bigr\}\le C\psi_{n,M}(s),
\end{eqnarray}
where
\[
\psi_{n,M}(s) = \frac{s}{n}\log\biggl(\frac{eM}{s} \biggr)
\wedge\sqrt{\frac{1}{n} \log\biggl(1+\frac{M}{\sqrt{n}} \biggr)}
\wedge1
\]
%
for $s\in\{1,\ldots, M\}$.


\end{theorem}

This theorem improves upon the rate of $s$-convex aggregation given in
Lounici \cite{Lounici07} by removing a redundant $(s/n)\log n$ term
present there. Note that \cite{Lounici07} considers the random design
regression model with Gaussian errors. Theorem~\ref{thmmainaggr} is
distribution-free and deals with bounded errors as all the results of
this paper; it can be readily extended to the case of sub-exponential
errors. By an easy modification of the minimax lower bound given in
\cite{Lounici07}, we get that $\psi_{n,M}(s)$ is the optimal rate for
the minimax regret on $\F_{\Theta^{\mathrm{C}}(s)}$ in our setting.
Analogous result for Gaussian regression with fixed design is proved in
\cite{rigollet2011exponential}.

\setcounter{rem}{3}
\begin{rem}\label{r4}
Inspection of the proofs shows that Theorems~\ref
{thmmainc}--\ref{thmmainaggr} as well as Theorem~\ref
{thapproxerror} below provide bounds on the risk not only in
expectation but also in deviation. For example, under the assumptions
of Theorem~\ref{thmmainVC}, along with (\ref{eq2thmmainVC}) we
obtain that there exists a constant $C>0$ depending only on $A$ such
that, for any $t>0$,
%
\begin{eqnarray}
\label{eq2thmmainVC-bis} \sup_{P_{XY}\in\cP} \Prob\biggl\{ \llVert
\tilde{f}-\eta
\rrVert^2 -\inf_{f\in\F} \llVert f-\eta\rrVert
^2 \ge C \biggl(\frac
{v}{n} \log\biggl(\frac{en}{v}
\biggr)+\frac{t}{n} \biggr) \biggr\}\le \mathrm{e}^{-t}.
\end{eqnarray}
The ``in deviation'' versions of Theorems~\ref{thmmainc}, \ref
{thmmainaggr} and~\ref{thapproxerror} are analogous and we skip
them for brevity. We also note that all the results trivially extend to
the case $\Y=[a,b]$, $\F\subseteq\{f\dvt a\le f \le b\}$, where
$-\infty<a<b<\infty$.
\end{rem}

\section{Adapting to approximation error rate of function class}\label{secapprox}


In Theorem~\ref{thmmainc}, we have shown that for $p>2$ our
estimator has the rate of $n^{-2/(2+p)}$ when $\eta\in\F$ and
achieves the rate of $n^{-1/p}$ if not. A natural question one can
ask is what happens if $\eta\notin\F$ but the approximation error
$\inf_{f \in\F} \llVert\eta- f\rrVert^2$ is small. This
can be viewed as
an intermediate setting between the pure statistical learning and
pure estimation. In such situation, one would expect to achieve
rates varying between $n^{-1/p}$ and $n^{-2/(2+p)}$ depending on how
small the approximation error is. This is indeed the case as
described in the next theorem.
%
\begin{theorem}\label{thapproxerror}
Let $\Y=[0,1]$, $\F\subseteq\{f\dvt 0\le f\le1\}$, and $\cH_2(\F
,\rho)\le A\rho^{-p}$, $\forall\rho>0$, for some constants
$A<\infty$, $p\geq2$. Consider\vspace*{1pt} an aggregation-of-leaders estimator
$\tilde{f}$ with the covering radius set as $\varepsilon= n^{-1/(2+p)}$. For this estimator and for any joint distribution $P_{XY}$
we have:
%
\begin{eqnarray}
\label{eqapproxerrorsmooth} \En\llVert\tilde{f} - \eta\rrVert^2 - \inf
_{f \in\F} \llVert f - \eta\rrVert^2 \le
C_p \bar\psi_{n,p}(\Delta),
\end{eqnarray}
where $\Delta^2 = \inf_{f \in\F} \llVert f-\eta\rrVert
^2$, $C_p>0$ is a
constant depending only on $p$ and $A$, and
%
\begin{eqnarray}
\label{defbarpsi} \bar\psi_{n,p}(\Delta)= \cases{ n^{-2/(2+p)},
&\quad
if $\Delta^2 \le n^{-2/(2+p)}$,
\vspace*{3pt}\cr
\Delta^2, &\quad
if $n^{-2/(2+p)} \le\Delta^2 \le n^{-1/p}$,
\vspace*{3pt}\cr
n^{-1/p}, &\quad if $\Delta^2 \ge n^{-1/p}$}
\end{eqnarray}
for $p>2$.
At $p=2$ the rate $\bar\psi_{n,p}(\Delta)$ is $n^{-1/2}\log n$
independently of $\Delta$.
\end{theorem}

The proof of this theorem is given in Section~\ref{secproofsthms}.

For particular cases $\Delta=0$ (well-specified model) or $\Delta=1$
(misspecified model), we recover the result of
Theorem~\ref{thmmainc} for $p>2$. Theorem~\ref{thapproxerror}
reveals that there is a smooth transition in terms of approximation
error rate in the intermediate regime between these two extremes.
Note also that the estimator ${\tilde f}$ in
Theorem~\ref{thapproxerror} is the same in all the cases; it is
defined by the aggregation-of-leaders procedure with $\varepsilon$
fixed as $n^{-1/(2+p)}$. Thus, the estimator is \emph{adaptive} to the
approximation error.

Theorem~\ref{thapproxerror} naturally suggests to study a minimax
problem which is more general than those considered in Statistical
Learning Theory or Nonparametric Estimation. Introduce the class of
$\Delta$-misspecified models
\[
\cP_\Delta(\F) = \Bigl\{P_{XY}\in\cP\dvt \inf
_{f \in\F} \llVert f - \eta\rrVert\le\Delta\Bigr\}, \qquad\Delta
\ge0, %
\]
and define the \textit{$\Delta$-misspecified regret} as
\[
V^\Delta_n(\F) = \inf_{\hat{f}} \sup
_{P_{XY}\in\cP_\Delta(\F
)} \Bigl\{ \En\llVert\hat{f}-\eta\rrVert^2 -
\inf_{f\in
\F} \llVert f-\eta\rrVert^2 \Bigr\} .
\]
Note that by definition, $V^\Delta_n(\F) = W_n(\F)$ when $\Delta=
0$ and $V^\Delta_n(\F) = V_n(\F)$ when $\Delta= 1$ (the diameter of
$\F$). In general, $V^\Delta_n(\F)$ measures the minimax regret when
we consider the statistical estimation problem with approximation
error at most $\Delta$. Theorem~\ref{thapproxerror} implies that
the rate of convergence of $\Delta$-misspecified regret admits the
bound $V^\Delta_n(\F)\le C_p\bar\psi_{n,p}(\Delta)$.

\section{Lower bounds}\label{seclower}

In this section, we show that the upper bounds obtained in
Theorems~\ref{thmmainc}, \ref{thmmainVC},
and~\ref{thapproxerror} cannot be improved. First, we exhibit a
VC-subgraph class $\F$ with VC-dimension at most $d$ such that
\[
W_n(\F)\geq C\frac{d}{n} \log\biggl(\frac{en}{d}
\biggr),
\]
where
$C>0$ is a numerical constant. In fact, we will prove a more
general lower bound, for the risk in probability rather than in
expectation.\vspace*{1pt}

In the next theorem, $\X=\{x^1,x^2,\ldots\}$ is a countable set of
elements and $\F$ is the following set of binary-valued functions on
$\X$:
\[
\F=\bigl\{f\dvt f(x) = a{\mathbf1} \{x\in W \}\mbox{ for some }
W\subset\X\mbox{ with }
\llvert W\rrvert\leq d\bigr\},
\]
where $a>0$, ${\mathbf1} \{\cdot\}$ denotes the indicator
function, $\llvert W\rrvert$ is
the cardinality of $W$, and $d$ is an integer. It is easy to check
that $\F$ is a VC-subgraph class with VC-dimension at most $d$.

%
\begin{theorem} \label{thlow1} Let $d$ be any integer such that $n\ge
d$, and $a=3/4$. Let the random pair $(X,Y)$ take values in $\X\times
\{0,1\}$. Then there exist a marginal distribution $\mu_X$ and
numerical constants $c,c'>0$ such that
\[
\inf_{\hat f} \sup_{\eta\in\F}{
\Prob}_\eta\biggl( \llVert\hat{f} - \eta\rrVert^2 \geq
c\frac{d}{n} \log\biggl(\frac
{en}{d} \biggr) \biggr)\geq
c',
\]
where ${\Prob}_\eta$ denotes the distribution of the $n$-sample $D_n$
when $\En(Y\mid X=x)=\eta(x)$.
\end{theorem}

The proof of Theorem~\ref{thlow1} is given in Section~\ref{secproofslower}.

We now exhibit a class $\F$ with polynomial growth of the empirical
entropy, for which the rates of minimax risk and minimax regret given
in Theorems~\ref{thmmainc} and~\ref{thapproxerror} cannot be
improved on any estimators. To state the result, we need some notation.
Let $\ell$ be the set of all real-valued sequences $(f_k, k=1,2,\ldots
)$. Denote by $\e_j$ the unit vectors in $\ell\dvt \e_j= ({\mathbf
1} \{k=j \},
k=1,2,\ldots)$, $j=1,2,\ldots.$ For $p>0$, consider the set
$B_p\triangleq\{f\in\ell\dvt \llvert f_j\rrvert\le
j^{-1/p}, j=1,2,\ldots\}$.

The next theorem provides lower bounds on $V_n(\F)$ and
$W_n(\F)$ when the $\varepsilon$-entropy of $\F$
behaves as $\varepsilon^{-p}$. It implies that the rates for $V_n(\F)$ and
$W_n(\F)$ in Theorem~\ref{thmmainc} are tight when $p< 2$. 

%
\begin{theorem}
\label{lemlowerrich}
Fix any $p>0$. Let $\F=\{f\in\ell\dvt f_j=(1+g_j)/2, \{g_j\}\in
B_p\}$ and let $\X=\{\e_1,\e_2,\ldots\}$ be the set of all unit
vectors in $\ell$. For any $\varepsilon>0$ 
we have
%
\begin{eqnarray}
\label{thlow21} \cH_2(\F,\varepsilon) \leq\biggl(\frac{A}{\varepsilon}
\biggr)^p,
\end{eqnarray}
where $A$ is a constant depending only on $p$. Furthermore, for this
$\F$, there exists an absolute positive constant $c$
such that the minimax risk satisfies, for any $n\ge1$,
%
\begin{eqnarray}
\label{thlow22} W_n(\F) \geq c n^{-2/(2+p)},
\end{eqnarray}
and the minimax regret satisfies, for any $p\geq2$ and any $n\ge1$,
%
\begin{eqnarray}
\label{thlow23} V_n(\F)& \geq& c n^{-1/(p-1)}.
\end{eqnarray}
%
\end{theorem}

The proof of Theorem~\ref{lemlowerrich} is given in Section~\ref
{secproofslower}. We remark that the lower bound (\ref{thlow23})
(for $p> 2$) holds, up to logarithmic factors, for \emph{any} class
satisfying the entropy growth $\Omega(\varepsilon^{-p})$, but we omit
the longer proof of this fact. We also remark that for $p>2$, the
$n^{-1/p}$ lower bound can be shown for any estimator taking values
within the class $\F$. Obtaining such a lower bound for any estimator
remains an open problem.

%

\section{Comparison with global ERM and with skeleton aggregation}\label{seccompare}

Among the methods of estimation designed to work under general
entropy assumptions on $\F$, the global ERM or the ERM on
$\varepsilon$-nets \cite{devroye87,BuescherKumar96,LugNob99} hold a
dominant place in the literature (see an
overview in Section~\ref{sechistorical}). Somewhat less studied
method is skeleton aggregation \cite{yang1999information}. In this
section, we discuss the deficiencies of these two previously known
methods that motivated us to introduce aggregation-of-leaders.

Recall that the aggregation-of-leaders procedure has three steps.
The first one is to find an empirical $\varepsilon$-net (that we will
call a skeleton) from the first subsample and partition the function
class based on the skeleton using the empirical distance on this
subsample. In the next step, using the second subsample we find
empirical risk minimizers within each cell of the partition.
Finally, we use the third sample to aggregate these ERM's. A simpler
and seemingly intuitive procedure that we will call the \textit
{skeleton aggregation} consists of steps one and three, but not two.
This method directly aggregates centers of the cells
$\hat{\F}_i^S(\varepsilon)$, that is, the elements ${\hat c}_i$ of the\vadjust{\goodbreak}
$\varepsilon$-net obtained from the first subsample $S$. Such kind of
procedure was studied by Yang and Barron
\cite{yang1999information} in the context of well-specified models.
The setting in \cite{yang1999information} is different from ours
since in that paper the $\varepsilon$-net is taken with respect to a
non-random metric and the bounds on the minimax risk $W_n(\F)$ are
obtained when the regression errors are Gaussian. Under this model,
\cite{yang1999information} provides the bounds not for skeleton
aggregation but for a more complex procedure that comprises an
additional projection in Hellinger metric. We argue that, while the
skeleton aggregation achieves the desired rates for well-specified
models (i.e., for the minimax risk), one cannot expect it to be
successful for the misspecified setting. This will explain why
aggregating ERM's in cells of the partition, and not simply
aggregating the centers of cells, is crucial for the success of the
aggregation-of-leaders procedure.

Let us first show why the skeleton aggregation yields the correct rates
for well-specified models (i.e., when $\eta\in\F$). Similarly
to~(\ref{mixture}), we define the skeleton aggregate ${\tilde
f}^{\mathrm{sk}} =\sum_{i=1}^N \theta_i{\hat c}_i $ 
as a sharp $\mathrm{MS}$-aggregate satisfying a bound analogous to (\ref
{eqaggregationbound}):
there exists a constant $C>0$ such that, for any $\delta>0$,
%
\begin{eqnarray}
\label{eqaggregationbound1} L\bigl({\tilde f}^{\mathrm{sk}}\bigr) \leq
\min
_{i=1,\ldots,N} L(\hat{c}_i) + C\frac
{\log(N/\delta) }{n}
\end{eqnarray}
with probability at least $1-\delta$ over the sample $S''$,
conditionally on $S$ (the subsample $S'$ is not used here). If the
model is well-specified, $L^*=L(\eta)$, and $\llVert f - \eta
\rrVert^2 = L(f)
-L^*$, $\forall f\in\F$. Hence, with probability $1-5\delta$,
%
\begin{eqnarray} \label{eqaggregationbound2}
\bigl\llVert{\tilde f}^{\mathrm{sk}} - \eta\bigr\rrVert^2 &=& L
\bigl({\tilde f}^{\mathrm{sk}}\bigr) -L^*\nonumber
\\[-1pt]
&\leq& \min_{i=1,\ldots,N} L(
\hat{c}_i)-L^* + C\frac{\log(N/\delta)
}{n}
\nonumber\\[-8pt]\\[-8pt]\nonumber
&=& \min_{i=1,\ldots,N} \llVert\hat{c}_i-\eta\rrVert
^2 + C\frac{\log(N/\delta
) }{n}
\\[-1pt]
&\le&2 \varepsilon^2 + C \biggl(\frac{\cH_2(\F,\varepsilon)}{n} +\frac
{\log(1/\delta) }{n} +
r^* + \beta\biggr)\nonumber
\end{eqnarray}
for $\beta=(\log(1/\delta)+\log\log n)/n$, and $r^*=r^*(\G)$ with
$\G=\{(f-g)^2\dvt f,g\in\F\}$, where we have used Lemma~\ref
{lemempbound} with $f=\hat{c}_i$, $f'=\eta$ and the fact that $\min
_{i=1,\ldots,N}d_S(\hat{c}_i,\eta)\le\varepsilon$ for any $\eta\in\F
$. The optimal choice of $\varepsilon$ in (\ref{eqaggregationbound2})
is given by the balance relation $n\varepsilon^2 \asymp\cH_2(\F
,\varepsilon)$ and it can be deduced from Lemma~\ref
{lemboundonr-star} that $r^* + \beta$ is negligible as compared to
the leading part $\mathrm{O}(\varepsilon^2+ \cH_2(\F,\varepsilon)/n)$ with this
optimal $\varepsilon$. In particular, we get from (\ref
{eqaggregationbound2}) combined with (\ref{eqrnstar-large}) and
(\ref{eqradpin02}), (\ref{radem}) that, under the assumptions of
Theorem~\ref{thmmainc}, $\sup_{\eta\in\F}\En\llVert{\tilde
f}^{\mathrm{sk}} - \eta\rrVert^2 \le
C n^{-2/(2+p)}$, $\forall p>0$.

Let us now consider the misspecified model setting (i.e., the
statistical learning framework). Here, the balance relation for the
skeleton aggregation takes the form $n\varepsilon\asymp\cH_2(\F
,\varepsilon)$, which yields suboptimal rates unless the class $\F$ is
finite. Indeed, without the assumption that the regression function
$\eta$ is in $\F$, we only obtain the bounds
%
\begin{eqnarray}
\label{eqaggregationbound4} L(\hat{c}_i)-L^* &=& \llVert\hat{c}_i -
\eta\rrVert^2 - \inf_{f \in\F} \llVert f - \eta\rrVert
^2
\nonumber\\[-9pt]\\[-9pt]\nonumber
& \le&2 \bigl( \llVert\hat{c}_i - \eta\rrVert- \llVert
\eta_{\F}- \eta\rrVert\bigr)+\frac{1}{n}\le2 \llVert
\hat{c}_i - \eta_{\F}\rrVert+\frac{1}{n},
\end{eqnarray}
where $\eta_{\F}\in\F$ is such that $\llVert\eta_{\F}- \eta
\rrVert^2
\le\inf_{f \in\F} \llVert f - \eta\rrVert^2+1/n$. The
crucial difference
from (\ref{eqaggregationbound2}) is that here $L(\hat{c}_i)-L^*$
behaves itself as a norm $ \llVert\hat{c}_i - \eta_{\F}\rrVert$ and
not as
a squared norm $ \llVert\hat{c}_i - \eta\rrVert^2$. Using
(\ref{eqaggregationbound4}) and arguing analogously to (\ref
{eqaggregationbound2}), we find that for misspecified models, with
probability $1-5\delta$,
%
\begin{eqnarray} \label{eqaggregationbound3}
L\bigl({\tilde f}^{\mathrm{sk}}\bigr) -L^* &\leq&2\min_{i=1,\ldots,N}
\llVert\hat{c}_i - \eta_{\F}\rrVert+ C\frac{\log(N/\delta)
}{n}\nonumber
\\
&\le&2\sqrt{2 \varepsilon^2 + C\bigl(r^* + \beta\bigr)} + C
\frac{\log(N/\delta
) }{n}
\\
&\le&2 \sqrt{2}\varepsilon+ C \biggl(\frac{\cH_2(\F,\varepsilon)}{n}
+\frac{\log(1/\delta) }{n} +
\sqrt{r^* + \beta} \biggr).\nonumber
\end{eqnarray}
Here, the optimal $\varepsilon$ is obtained from the tradeoff of
$\varepsilon$ with $\cH_2(\F,\varepsilon)/n$. As a result, we only get
the suboptimal rate $n^{-1/(p+1)} + \mathrm{O}(\sqrt{r^* + \beta})$ for the
excess risk of ${\tilde f}^{\mathrm{sk}}$ under the assumptions of
Theorem~\ref{thmmainc}. While the above argument is based on upper
bounds, it is possible to construct a simple scenario where $\eta,
\eta_\F$ and some $\hat{c}_i$ are on a line, $\llVert\eta_\F
-\hat
{c}_i\rrVert=\mathrm{O}(\varepsilon)$, and no other element $\hat{c}_j$ is
closer to
$\eta$ than $\hat{c}_i$. For such a setup, $L(\hat{c}_i)-L^*$ is of
the order of $\varepsilon$ and no convex combination of $\hat{c}_j$ can
improve upon $\hat{c}_i$. This indicates that introducing least
squares estimators over cells of the partition (the second step of our
procedure) is crucial in getting the right rates.

%
\begin{table}[b]
\tabcolsep=0pt
\caption{Summary of rates for misspecified case}\label{table1}
\begin{tabular*}{\tablewidth}{@{\extracolsep{\fill}}@{}llll@{}}
\hline
 & Aggregation- &  & \\
Regime & of-leaders & Skeleton aggregation & ERM\\
\hline
Finite: $\llvert\F\rrvert = M$ & $\frac{\log M}{n}$ & $\frac{\log M}{n}$ & $\sqrt{\frac{\log M}{n}}$
\\[6pt]
Parametric: $\operatorname{VC}(\F) =v \le n$ & $\frac{v \log(en/v)}{n}$& $\sqrt {\frac{v \log(en/v)}{n}}$ & $\sqrt{\frac{v}{n}}$
\\[6pt]
Nonparametric: $\cH_2(\F,\varepsilon) = \varepsilon^{-p}, $ & & & \\
\quad $p \in(0,2)$ & $n^{-2/(2+p)}$ & $n^{-1/(p+1)}\vee n^{-1/2}(\log n)^{3/2}$ & $n^{-1/2}$ \\
\quad $p \in(2,\infty)$ & $n^{-1/p}$ & $n^{-1/(p+1)}$ &$n^{-1/p}$ \\
\hline
\end{tabular*}
\end{table}

We can now compare the following three estimators. First, we consider
the global ERM over $\F$ defined~by
%
\begin{equation}
\label{deferm} \hat{f}^{\mathrm{erm}} \in\mathop{\operatorname{argmin}}_{f\in\F}
\frac{1}{n}\sum_{(x,y)\in
S'} \bigl(f(x)-y
\bigr)^2,
\end{equation}
second -- the skeleton aggregate ${\tilde f}^{\mathrm{sk}}$ and, finally,
the proposed aggregation-of-leaders estimator ${\tilde f}$. Table~\ref{table1}
summarizes the convergence rates of the expected excess risk $\cE_{\F
}(\hat{f})$ for $\hat{f}\in\{{\tilde f}, {\tilde f}^{\mathrm{sk}},\hat
{f}^{\mathrm{erm}}\}$ in misspecified model setting, that is, upper bounds
on the minimax regret.

The rates for finite $\F$ in Table~\ref{table1} are obtained in a trivial way by
taking the skeleton that coincides with the $M$ functions in the class
$\F$. In parametric and nonparametric regime, the rates for the
proposed method are taken from Theorems~\ref{thmmainc} and~\ref
{thmmainVC}, while for the skeleton aggregate they follow from (\ref
{eqaggregationbound3}) with optimized $\varepsilon$ combined with the
bounds on $r^*$ in Lemma~\ref{lemboundonr-star} and in (\ref
{eqradpin02}), (\ref{radem}) below. The rate $\sqrt{v/n}$ for the
excess risk of ERM in parametric case is well-known, cf., for example,
\cite{bousquet2002concentration,bartlett2006notes}. For the
nonparametric regime, the rates for ERM in Table~\ref{table1} follow from
Lemma~\ref{lemerm} and the bounds on $\Rad_n(\F)$ in (\ref
{eqradpin02}) and (\ref{radem}) below. Moreover, for finite $\F$,
it can be shown that the slow rate $\sqrt{\frac{\log M}{n}}$ cannot
be improved neither for ERM, nor for any other selector, that is, any
estimator with values in $\F$, cf. \cite{jrt-lma-08}.

In conclusion, for finite class $\F$ aggregation-of-leaders and
skeleton aggregation achieve the excess risk rate $\frac{\log M}{n}$,
which is known to be optimal \cite{tsybakov03optimal}, whereas the
global ERM has a suboptimal rate. For a very massive class $\F$, when
the empirical entropy grows polynomially as $\varepsilon^{-p}$ with $p
\ge2$ both ERM and aggregation-of-leaders enjoy similar guarantees of
rates of order $n^{ - 1/p}$ while the skeleton aggregation only gets a
suboptimal rate of $n^{-1/(p+1)}$. For all other cases, while
aggregation-of-leaders is optimal, both ERM and skeleton aggregation
are suboptimal. Thus, in the misspecified case, skeleton aggregation is
good only for very meager (finite) classes while ERM enjoys optimality
only for the other extreme -- massive nonparametric classes. Note also
that, unless $\F$ is finite, skeleton aggregation does not improve
upon ERM in the misspecified case.

Turning to the well-specified case, both aggregation-of-leaders and
skeleton aggregation achieve the optimal rate 
for the minimax risk while the global ERM is, in general, suboptimal.

\section{Historical remarks and comparison with previous work}\label{sechistorical}


The role of entropy and capacity \cite{kolmogorov1959} in establishing
rates of estimation has been recognized for a long time, since the work
of Le Cam~\cite{lecam1973convergence}, Ibragimov and Has'minski{\u\i}
\cite{ibragimov1980estimate} and Birg\'e \cite{Birge83}. This was
also emphasized by Devroye~\cite{devroye87} and Devroye \textit{et al}.~\cite
{DevrGyorLug96} in the study ERM on $\varepsilon$-nets. The common point
is that optimal rate is obtained as a solution to the balance equation
$n\varepsilon^2 = \cH(\varepsilon)$, with an appropriately chosen
non-random entropy $\cH(\cdot)$. 
Yang and Barron \cite{yang1999information} present a general approach
to obtain lower bounds from global (rather than local) capacity
properties of the parameter set. Once again, the optimal rate is shown
to be a solution to the bias-variance balance equation described above,
with a generic notion of a metric on the parameter space and non-random
entropy. Under the assumption that the regression errors are Gaussian,
\cite{yang1999information} also provides an achievability result via a
skeleton aggregation procedure complemented by a Hellinger projection
step. 
Van de Geer~\cite{Gee90} invokes the empirical entropy rather than the
non-random entropy to derive rates of estimation in regression problems.

In all these studies, it is assumed that the unknown density,
regression function, or parameter belongs to the given class, that is,
the model is well-specified. In parallel to these developments, a line
of work on pattern recognition that can be traced back to Aizerman,
Braverman and Rozonoer \cite{aizerman1970method} and Vapnik and
Chervonenkis \cite{vapnik1974theory} focused on a different objective,
which is characteristic for Statistical Learning. Without assuming a
form of the distribution that encodes the relationship between the
predictors and outputs, the goal is formulated as that of performing as
well as the best within a given set of rules, with the excess risk as
the measure of performance (rather than distance to the true underlying
function). Thus, no assumption is placed on the underlying
distribution. In this form, the problem can be cast as a special case
of stochastic optimization and can be solved either via recurrent
(e.g., gradient descent) methods or via empirical risk minimization.
The latter approach
leads to the question of uniform convergence of averages to
expectations, also called the uniform Glivenko--Cantelli property. This
property is, once again, closely related to entropy of the class, and
sufficient conditions have been extensively studied (see \cite
{Dudley78,Pollard84,Dudley99,Dudley87,DudGinZin91} and references therein).

For parametric classes with a polynomial growth of covering numbers,
uniform convergence of averages to expectations with the $\sqrt{(\log
n)/n}$ rate has been proved by Vapnik and Chervonenkis \cite
{vc-ucrfep-71,VapChe68,vapnik1974theory}. In the context of
classification, they also obtained a faster rate showing $\mathrm{O}((\log n)
/n)$ convergence when the minimal risk $L^*=0$. For regression
problems, similar fast rate when \mbox{$L^*=0$} can be shown (it can be
deduced after some argument from Assertion 2 on page 204 in~\cite
{Vap82}; an exact formulation is available, e.g., in \cite
{srebro2010smoothness}). Lee, Bartlett and Williamson \cite
{lbw-iclsl-98} showed that the excess risk of the least squares
estimator on $\F$ can attain the rate $\mathrm{O}((\log n)/n)$ without the
assumption $L^*=0$. Instead, they assumed that the class $\F$ is
convex and has finite pseudo-dimension. Additionally, it was shown that
the $n^{-1/2}$ rate cannot be improved if the class is non-convex and
the estimator is a selector (i.e., forced to take values in $\F$).
In particular, the excess risk of ERM and of any selector on a finite
class $\F$ cannot decrease faster than $\sqrt{(\log\llvert\F
\rrvert)/n}$ \cite
{jrt-lma-08}. Optimality of ERM for certain problems is still an open question.

Independently of this work on the excess risk in the distribution-free
setting of statistical learning, Nemirovskii \cite
{nemirovski2000topics} proposed to study the problem of aggregation, or
mimicking the best function in the given class, for regression models.
Nemirovskii~\cite{nemirovski2000topics} outlined three problems: model
selection, convex aggregation, and linear aggregation. The notion of
optimal rates of aggregation based on the minimax regret is introduced
in \cite{tsybakov03optimal}, along with the derivation of the optimal
rates for the three problems. In the following decade, much work has
been done on understanding these and related aggregation problems \cite
{Yang04,jntv-raemdaa-05,jrt-lma-08,Lounici07,rigollet2011exponential}.
For recent developments and a survey we refer to \cite{LecueHabil,RigTsySTS12}.

In parallel with this research, the study of the excess risk blossomed
with the introduction of Rademacher and local Rademacher complexities
\cite
{Kol01,KolPan2000,bm-rgcrbsr-02,bousquet2002some,bbm-lrc-05,koltchinskii2006local}.
These techniques provided a good understanding of the behavior of the
ERM method. In particular, if $\F$ is a \emph{convex} subset of
$d$-dimensional space, Koltchinskii \cite
{koltchinskii2006local,koltchinskii2011oracle} obtained a sharp oracle
inequality with the correct rate $d/n$ for the excess risk of least
squares estimator on $\F$. Also, for convex $\F$ and $p\in(0,2)$,
the least squares estimator on $\F$ attains the correct excess risk
rate $n^{-2/(p+2)}$ under the assumptions of Theorem~\ref
{thmmainc}. This can be deduced from Theorem~5.1 in \cite
{koltchinskii2011oracle}, remarks after it and in Example~4 on page 87
of \cite{koltchinskii2011oracle}. However, the convexity assumption
appears to be crucial; without this assumption Koltchinskii \cite
{koltchinskii2011oracle}, Theorem 5.2, obtains for the least squares
estimator only a non-sharp inequality with leading constant $C>1$, cf.
(\ref{eqnonexactoracle}). As follows from the results in
Section~\ref{secmain} our procedure overcomes this problem.

Among a few of the estimators considered in the literature for general
classes $\F$, empirical risk minimization on $\F$ has been one of the
most studied. As mentioned above, ERM and other selector methods are
suboptimal when the class $\F$ is finite. 
For the regression setting with finite~$\F$, the approach that was
found to achieve the optimal rate for the excess risk in expectation is
through exponential weights with averaging of the trajectory~\cite
{yang1999information,catoni,jrt-lma-08,DalTsy12}. However, Audibert
\cite{audibert2007progressive} showed that, for the regression with
random design, exponential weighting is suboptimal when the error is
measured by the probability of deviation rather than by the expected
risk. He proposed an alternative method, optimal both in probability
and in deviation, which involves finding an ERM on a star connecting a
global ERM and the other $\llvert\F\rrvert-1$ functions.
In \cite{lecue2009aggregation}, the authors exhibited another
deviation optimal method which involves sample splitting. The first
part of the sample is used to localize a convex subset around ERM and
the second -- to find an ERM within this subset. 
Recently yet another procedure achieving the deviation optimality has
been proposed in \cite{LecueRigollet}. It is based on a penalized
version of exponential weighting and extends the method of \cite
{dai2012deviation} originally proposed for regression with fixed design.
The methods of \cite
{audibert2007progressive,lecue2009aggregation,LecueRigollet} provide
examples of sharp MS-aggregates that can be used at the third step of
our procedure.

We close this short summary with a connection to a different
literature. In the context of prediction of deterministic individual
sequences with logarithmic loss, Cesa-Bianchi and Lugosi \cite
{cbl-wcbllp-01} considered regret with respect to rich classes of
``experts''. They showed that mixture of densities is suboptimal and
proposed a two-level method where the rich set of distributions is
divided into small balls, the optimal algorithm is run on each of these
balls, and then the overall output is an aggregate of outputs on the
balls. They derived a bound where the upper limit of the Dudley
integral is the radius of the balls. This method served as an
inspiration for the present work.

\section{Proofs of Theorems~\texorpdfstring{\protect\ref{thmmainc}}{2}--\texorpdfstring{\protect\ref{thmmainaggr}}{4} 
and \texorpdfstring{\protect\ref{thapproxerror}}{5}}
\label{secproofsthms}


We first state some auxiliary lemmas.

%
\begin{lemma}
\label{lemboundonr-star}
The following values can be taken as localization radii $r^*=r^*(\G)$
for $\G=\{(f-g)^2\dvt f,g\in\F\}$.
\begin{longlist}[(iii)]
\item[(i)] For any class $\F\subseteq\{f\dvt 0\le f\le1\}$, and $n\ge2$,
%
\begin{eqnarray}
\label{eqrnstar-large} r^* = C\log^3(n) \Rad_n^2(
\F).
\end{eqnarray}
\item[(ii)] If $\F\subseteq\{f\dvt 0\le f\le1\}$ and the empirical
covering numbers exhibit polynomial growth
$\sup_{S\in\Z^n}\cN_2(\F,\rho,S) \leq(\frac{A}{\rho
} )^v$
for some constants $A<\infty$, $v>0$, then
\[
r^*= C\frac{v}{n}\log\biggl(\frac{en}{v} \biggr)
\]
whenever $n\ge C_Av$ with $C_A>1$ large enough depending only on $A$.
\item[(iii)] If $\F$ is a finite class with $\llvert\F\rrvert\ge2$,
\[
r^*= C\frac{\log\llvert\F\rrvert}{n} .
\]
\end{longlist}
\end{lemma}
The proof of this lemma is given in the \hyperref[secappendix]{Appendix}. The following lemma
is a direct consequence of Theorem~\ref{thmlocalization} proved in
the \hyperref[secappendix]{Appendix}.
%
\begin{lemma}
\label{lemempbound} For any class $\F\subseteq\{f\dvt 0\le f\le1\}$
and $\delta>0$, with probability at least $1-4\delta$,
%
\begin{eqnarray}
\label{eqlemempbound} \bigl\llVert f-f'\bigr\rrVert^2
\leq2d_S^2\bigl(f,f'\bigr) + C\bigl(r^* +
\beta\bigr)\qquad\forall f,f'\in\F,
\end{eqnarray}
where $\beta=(\log(1/\delta)+\log\log n)/n$, and $r^*=r^*(\G)$ for
$\G=\{(f-g)^2\dvt f,g\in\F\}$.
\end{lemma}

We will also use the following bound on the Rademacher average in terms
of the empirical entropy \cite{bartlett2006notes,srebro2010smoothness}.
%
\begin{lemma}
\label{lemradent} For any class $\F\subseteq\{f\dvt 0\le f\le1\}$,
%
\begin{eqnarray}
\label{radent} \hat{\Rad}_n(\F,S) &\leq&\inf_{\alpha\geq0}
\biggl\{ 4\alpha+ \frac{12}{\sqrt{n}}\int_{\alpha}^1
\sqrt{\log\cN_2(\F,\rho,S)} \,\mathrm{d}\rho\biggr\} .
\end{eqnarray}
\end{lemma}

%
\begin{pf*}{Proof of Theorem~\ref{thmmainc}}
Consider the case $p\in(0,2)$. 
Assume without loss of generality that $A=1$, that is, $\sup_{S\in\Z
^n}\log\cN_2(\F, \rho,S) \le\rho^{-p}$.
For $p\in(0,2)$, the bound (\ref{radent}) with $\alpha=0$ combined
with (\ref{eqrnstar-large}) yields
%
\begin{eqnarray}
\label{eqradpin02} {\Rad}_n(\F) &\leq&\frac{12}{\sqrt{n}(1-p/2)} ,
\qquad r^* \leq
C \frac{(\log n)^3}{n}
\end{eqnarray}
for some absolute constant $C$. Thus,
%
\begin{eqnarray}
\label{thm2proof1} \gamma&\le& C \biggl(\varepsilon+ \frac{(\log
n)^{3/2}}{\sqrt{n}} + \sqrt{
\frac{\log(1/\delta)}{n}} \biggr) ,
\\
\gamma\sqrt{r^*} &\leq& C(\log n)^{3/2} \biggl( \frac{\varepsilon
}{\sqrt{n}}+
\frac{(\log n)^{3/2}}{n} + \frac{\sqrt{\log(1/\delta
)}}{n} \biggr) .
\end{eqnarray}
These inequalities together with (\ref{1}) and (\ref{2}) yield
that for $0<\delta<1/2$, with probability at least $1-2\delta$,
%
\begin{equation}
\label{thm2proof2} L (\tilde{f} )-L^* \leq C \biggl(\frac{\varepsilon
^{-p}}{n} +
\frac{\log(1/\delta) }{n} + \gamma\sqrt{r^*} + \frac{\gamma
^{1-p/2}}{\sqrt{n}} \biggr) .
\end{equation}
The value of $\varepsilon$ minimizing the right-hand side in (\ref
{thm2proof2}) is $\varepsilon= n^{-1/(2+p)}$, which justifies the choice
made in the theorem.
Notably, the logarithmic factor arising from $r^*$ only appears
together with the lower order terms and the summand $\gamma\sqrt
{r^*}$ does not affect the rate.
For $\varepsilon= n^{-1/(2+p)}$ the right-hand side of (\ref
{thm2proof2}) is bounded by $Cn^{-2/(2+p)}$ ignoring the terms
with $\log(1/\delta)$ that disappear when passing from the bound in
probability to that in expectation. Thus, the expected excess risk is
bounded by $Cn^{-2/(2+p)}$, which proves (\ref{4}) for $p\in(0,2)$.

Next, consider the case $p>2$. From (\ref{radent}) with $\alpha
=n^{-1/p}$, $\hat{\Rad}_n(\F,S)\leq Cn^{-1/p}$ and $r^* = (\log n)^3
n^{-2/p}$. Choosing $\varepsilon=n^{-1/(2+p)}$,
\[
\gamma\sqrt{r^*} \leq C\bigl(\varepsilon\sqrt{r^*}+r^*+\sqrt{\beta
r^*}\bigr)\leq C
n^{-1/p} .
\]
The first statement of the theorem follows from (\ref{2}) with the
choice $\alpha=n^{-1/p}$ and by noting that $\frac{\varepsilon
^{-p}}{n}$ is of the lower order than $n^{-1/p}$. The case of $p=2$
follows similarly (see proof of Theorem~\ref{thapproxerror}). The
second part of the theorem follows from Theorem~\ref{thapproxerror}.
\end{pf*}

%
\begin{pf*}{Proof of Theorem~\ref{thmmainVC}}
Throughout this proof, $C$ is a generic notation for positive constants
that may depend only on $A$.
Since $\varepsilon=n^{-1/2}$ the expression for $r^*$ in Lemma~\ref
{lemboundonr-star}(ii) leads to the bounds $\gamma\leq C
(\sqrt{\frac{v\log(en/v)}{n}}+ \sqrt{ \frac{\log(1/\delta
)}{n}} )$, and $\gamma\sqrt{r^*} \leq C (\frac{v\log
(en/v)}{n}+ \frac{\log(1/\delta)}{n} )$. Next, since $ \cN
_2(\F,\rho,S')\le\max\{1, (A/\rho)^v\}$ we get
\begin{eqnarray*}
\frac{1}{\sqrt{n}}\int_{0}^{C\gamma} \sqrt{\log
\cN_2\bigl(\F, \rho, S'\bigr)} \,\mathrm{d}\rho&\leq&
\sqrt{\frac{v}{n}} \int_{0}^{C\gamma/A\wedge1} \sqrt{\log
(1/t)}\,\mathrm{d}t
\\
&\le& C \sqrt{\frac{v}{n}} \gamma\sqrt{\log({C}/{\gamma} )\vee1},
\end{eqnarray*}
where the last inequality is due to~(\ref{eqvcintegral}). We assume
w.l.o.g. that in the last expression $C$ is large enough to guarantee
that the function $\gamma\mapsto\gamma\sqrt{\log({C}/{\gamma
} )\vee1}$ is increasing, so that we can replace $\gamma$ by
the previous upper bound. This yields, after some algebra,
\[
\gamma\sqrt{\log({C}/{\gamma} )\vee1}\le C \biggl(\frac{\sqrt{v}\log
(en/v)}{\sqrt{n}} +
\frac{\sqrt{\log(1/\delta
)\log(en/v)}}{\sqrt{n}} \biggr) %
\]
if $n\ge Cv$ for $C$ large enough. The above inequalities together with
(\ref{1}) and (\ref{2}) imply
that, with probability at least $1-2\delta$,
\begin{eqnarray*}
L(\tilde{f})-L^* \leq C \biggl(\frac{\log\cN_2(\F
,\varepsilon,S)}{n} + \frac{v\log(en/v)}{n} +
\frac{\log(1/\delta)}{n} \biggr) .
\end{eqnarray*}
Using that $\log\cN_2(\F,\varepsilon,S)\le\max\{1,(A/\varepsilon)^v\}
$ and integrating over $\delta$ we get the desired bound for the
expected excess risk $\En L(\tilde{f})-L^*$.
\end{pf*}

\begin{pf*}{Proof of Theorem~\ref{thmmainaggr}}
By definition of the estimator, for any fixed integer $m\le s$ and $\nu
$ such that $\llvert\nu\rrvert=m$ we first construct the
least squares estimators
over the cells $\F_{\nu,m}$:
%
\begin{equation}
\label{defest1} \hat{f}_{\nu,m}^{S,S'} \in\mathop{\operatorname{argmin}}
_{f\in\F_{\nu,m}} \frac
{1}{n}\sum_{(x,y)\in S'}
\bigl(f(x)-y\bigr)^2.
\end{equation}
Since $\F_{\nu,m}$ is a convex hull of $m$ functions we can apply
\cite{Lecue2013} to get that for any $t>0$, with probability at least
$1-\mathrm{e}^{-t}$,
%
\begin{eqnarray}
\label{eqproofaggr1} L \bigl(\hat{f}_{\nu,m}^{S,S'} \bigr) \le\inf
_{f\in\F_{\nu,m}} L(f) + C(\tilde\psi_{m,n} + t/n),
\end{eqnarray}
where
\[
\tilde\psi_{m,n} \triangleq\frac{m}{n} \wedge\sqrt{
\frac{1}{n}\log\biggl(1+\frac{m}{\sqrt{n}} \biggr)} . %
\]
Thus, the event $E$ where (\ref{eqproofaggr1}) holds simultaneously
for all $(m,\nu)\in{\mathcal I}=\{(m,\nu)\dvt m=1,\ldots,s, \llvert\nu
\rrvert=m\}
$ is of probability at least $1-N\mathrm{e}^{-t}$.
Here, $N=\llvert{\mathcal I}\rrvert$. Choose now $t=\log
(N/\delta)$. Then, on the
intersection of $E$ with the event where (\ref{eqaggregationbound})
holds we have that, with probability at least $1-2\delta$,
%
\begin{eqnarray}
\label{eqproofaggr2} L (\tilde{f} ) &\le&\inf_{f\in\F_{\Theta^{\mathrm
{C}}(s)}} L(f) + C \biggl(
\tilde\psi_{s,n} + \frac{\log(N/\delta)}{n} \biggr)
\nonumber\\[-8pt]\\[-8pt]\nonumber
& \le&\inf_{f\in\F_{\Theta^{\mathrm{C}}(s)}} L(f)+ C \biggl(
\frac{s}{n}\log\biggl(\frac{eM}{s} \biggr)+ \frac{\log(1/\delta
)}{n} \biggr),
\nonumber
\end{eqnarray}
where we have used the inequalities $\tilde\psi_{m,n} \le\tilde\psi
_{s,n}$, $\forall m\le s$, and $N=\sum_{m=1}^s {M\choose m} \leq
(\frac{eM}{s} )^s$. On the other hand, for the least
squares estimator ${\hat f}^{\mathrm{C}}$ on the convex hull of all
$f_1,\ldots,f_M$, using again the result of \cite{Lecue2013} we have
that for any $u>0$, with probability at least $1-\mathrm{e}^{-u}$,
%
\begin{eqnarray}
\label{eqproofaggr3} L \bigl({\hat f}^{\mathrm{C}} \bigr) &\le&\inf
_{f\in\F_{\Theta^{\mathrm{C}}} } L(f) + C(\tilde\psi_{M,n} + u/n)
\nonumber\\[-8pt]\\[-8pt]\nonumber
&\le&\inf_{f\in\F_{\Theta^{\mathrm{C}}(s)} } L(f) + C \biggl(\sqrt{\frac
{1}{n}\log
\biggl(1+\frac{M}{\sqrt{n}} \biggr)}+ \frac
{u}{n} \biggr).
\end{eqnarray}
Now, we aggregate only two estimators, $\tilde{f}$ and ${\hat
f}^{\mathrm{C}}$ to obtain the final aggregate $\tilde{f}^*$. This
yields, in view
of (\ref{eqaggregationbound}) with $N=2$, (\ref{eqproofaggr2}),
and (\ref{eqproofaggr3}) with $u=\log(1/\delta)$, that with
probability at least $1-4\delta$,
\begin{eqnarray*}
L \bigl(\tilde{f}^* \bigr) &\le&\min\bigl\{L (\tilde{f} ), L \bigl({\hat
f}^{\mathrm{C}} \bigr) \bigr\} + C \frac
{\log(2/\delta)}{n}
\\
&\le&\inf_{f\in\F_{\Theta^{\mathrm{C}}(s)} } L(f)
+C \biggl(\min\biggl\{\frac{s}{n}\log\biggl(\frac{eM}{s} \biggr),
\sqrt{\frac{1}{n}\log\biggl(1+\frac{M}{\sqrt{n}} \biggr)} \biggr\} +
\frac{\log(1/\delta)}{n} \biggr) ,
\end{eqnarray*}
which immediately implies the desired bound for the expected excess
risk $\En L(\tilde{f}^*)-\break \inf_{f\in\F_{\Theta^{\mathrm{C}}(s)} }
L(f)$.
\end{pf*}


\begin{pf*}{Proof of Theorem~\ref{thapproxerror}}
Without loss of generality assume in this proof that $A=1$, that is,
that $\sup_{S\in\Z^n}\log\cN_2(\F, \rho,S) \le\rho^{-p}$.
Using (\ref{radent}) we bound ${\Rad}_n(\F)$ for $p>2$ as follows:
\begin{eqnarray}
\nonumber
{\Rad}_n(\F) \le\inf_{\alpha\geq0} \biggl\{ 4
\alpha+ \frac
{12}{\sqrt{n}} \int_{\alpha}^1
\rho^{-p/2}\,\mathrm{d}\rho\biggr\} &\leq&\inf_{\alpha\geq0} \biggl\{ 4
\alpha+ \frac{24}{\sqrt
{n}(p-2)}\alpha^{-(p-2)/2} \biggr\} .
\end{eqnarray}
For $p>2$, the balance equation $\alpha=n^{-1/2}\alpha^{-(p-2)/2}$
yields $\alpha=n^{-1/p}$. This and (\ref{eqrnstar-large}) lead to
the bounds
%
\begin{eqnarray}
\label{radem} {\Rad}_n(\F) & \le& Cn^{-1/p}, \qquad r^* \leq
C (\log n)^3 n^{-2/p}.
\end{eqnarray}
For $p=2$, choosing $\alpha= n^{-1/2}$,
%
\begin{eqnarray}
\label{radem-p-2} {\Rad}_n(\F) & \le& Cn^{-1/2}\log n, \qquad r^*
\leq C (\log n)^5 n^{-2/p}.
\end{eqnarray}

Consider the case $p>2$. Let $\eta_\F\in\F$ be such that $\llVert\eta
_\F- \eta\rrVert^2 \le\inf_{f \in\F} \llVert
f - \eta\rrVert^2+1/n$.
Lemma~\ref{lemempbound}, (\ref{eqrnstar-large}) and (\ref{radem})
imply that, with probability at least $1-4\delta$, for all $i=1,\ldots,N$,
\begin{eqnarray*}
\nonumber
\bigl\llVert\hat{f}_i^{S,S'} - \eta\bigr\rrVert
^2 - \inf_{f\in\F
}\llVert f - \eta\rrVert
^2 &\le&2 \bigl\llVert\hat{f}_i^{S,S'} -
\eta_\F\bigr\rrVert^2 + \llVert\eta_\F-
\eta\rrVert^2 +1/n
\\
& \le&4 d_S^2\bigl(\hat{f}_i^{S,S'},
\eta_\F\bigr) + \llVert\eta_\F- \eta\rrVert
^2 + C\bigl(r^*+\beta\bigr)+1/n
\nonumber
\\
& \le&4 d_S^2\bigl(\hat{f}_i^{S,S'},
\eta_\F\bigr) + \Delta^2 + C \biggl(\frac{(\log n)^3}{n^{2/p}} +
\frac{\log(1/\delta)}{n} \biggr) .%
\end{eqnarray*}
Since $\min_{i=1,\ldots,N}d_S(\hat{f}_i^{S,S'}, \eta_\F)\le
2\varepsilon$ and $\varepsilon= n^{-1/(2+p)}$ we get that, with probability
at least $1-4\delta$,
%
\begin{eqnarray}\label{thapprox1}
\min_{i=1,\ldots,N}L\bigl(\hat{f}_i^{S,S'}
\bigr) - L^*& =& \min_{i=1,\ldots,N}\bigl\llVert\hat{f}_i^{S,S'}
- \eta\bigr\rrVert^2 - \inf_{f\in
\F}\llVert f - \eta
\rrVert^2
\nonumber\\[-8pt]\\[-8pt]\nonumber
& \le&\Delta^2 + C \biggl( n^{-2/(2+p)} + \frac{\log
(1/\delta)}{n}
\biggr).
\end{eqnarray}
Further, Lemma~\ref{lemerm} and (\ref{radem}) imply that, with
probability at least $1-2\delta$,
\begin{eqnarray*}
\bigl\llVert\hat{f}_i^{S,S'} - \eta\bigr\rrVert
^2 - \inf_{f\in\hat
{\F}_i^S}\llVert f - \eta\rrVert
^2 &\le& C \biggl(\Rad_n\bigl(\hat{\F}_i^S
\bigr) + \frac{\log
(1/\delta)}{n} \biggr)
\\
&\le& C \biggl(n^{-1/p} + \frac{\log(1/\delta)}{n} \biggr) .
\end{eqnarray*}
Combining this bound with (\ref{thapprox1}) we can conclude that,
with probability at least $1-6\delta$,
\begin{eqnarray*}
\min_{i=1,\ldots,N}L\bigl(\hat{f}_i^{S,S'}\bigr)
- L^*& =& \min_{i=1,\ldots,N}\bigl\llVert\hat{f}_i^{S,S'}
- \eta\bigr\rrVert^2 - \inf_{f\in
\F}\llVert f - \eta
\rrVert^2
\\
&\le& C \biggl(\min\bigl( n^{-2/(2+p)} + \Delta^2 ,
n^{-1/p} \bigr)+\frac{\log(1/\delta)}{n} \biggr) .
\end{eqnarray*}
Together with (\ref{eqaggregationbound}), this yields the next bound
that holds with probability at least $1-7\delta$:
\begin{eqnarray*}
\llVert\tilde{f} - \eta\rrVert^2 - \inf_{f \in\F}
\llVert f - \eta\rrVert^2 &=& L (\tilde{f} ) -L^*
\\
&\le& C \biggl(\frac{A\varepsilon^{-p}}{n} +\min\bigl(n^{-2/(2+p)} +
\Delta^2 , n^{-1/p} \bigr) +\frac{\log(1/\delta
)}{n} \biggr)
\\
& \le& C \biggl( n^{-2/(2+p)} + \min\bigl( n^{-2/(2+p)} +
\Delta^2 , n^{-1/p} \bigr)+\frac{\log(1/\delta)}{n} \biggr) ,
\end{eqnarray*}
and~(\ref{eqapproxerrorsmooth}) follows. For $p=2$, the above bound
gains a factor $\log n$ in front of $n^{-1/p}$ only.
\end{pf*}

\section{Proof of Theorem~\texorpdfstring{\protect\ref{lemmain}}{1}}\label{secanalysis}

We start with the following bound on the risk of least squares
estimators in terms of Rademacher complexity.
%
\begin{lemma}
\label{lemerm}
Let 
$\F$ be a class of measurable functions from $\X$ to $[0,1]$. Then,
for any $t>0$, with probability at least $1-2\mathrm{e}^{-t}$, the least squares
estimator $\hat{f}^{\mathrm{erm}}$ on $\F$ based on a sample $S'$ of size~$n$ (cf. (\ref{deferm})) satisfies
\begin{eqnarray*}
L\bigl(\hat{f}^{\mathrm{erm}}\bigr) \le L^* + C \hat{\Rad}_n\bigl(
\ell\circ\F, S'\bigr) + \frac{C t}{n} .
\end{eqnarray*}
%
\end{lemma}

The proof of this lemma is given in the \hyperref[secappendix]{Appendix} and is based on
combination of results from~\cite{bbm-lrc-05}. Note that here we have
both the remainder term of the order $1/n$ and the leading constant 1,
which is crucial for our purposes. 

Using Lemma~\ref{lemerm} with $\F=\hat{\F}_i^S$ and the union
bound, we obtain that, with probability at least $1-2N\mathrm{e}^{-t}$, for all
$i=1,\ldots,N$,
%
\begin{eqnarray}
\label{eqfirstbound} L \bigl(\hat{f}_i^{S,S'} \bigr) \leq\inf
_{f\in\hat{\F}_i^S} L(f) + C\hat{\Rad}_n\bigl(\ell\circ\hat{
\F}_i^S,S'\bigr) + C t/n.
\end{eqnarray}
Recall that $N=\cN_2(\F,\varepsilon, S)$. Setting $t=\log(4N/\delta)$
and using (\ref{eqfirstbound}) and (\ref{eqaggregationbound}) we
obtain that, with probability at least $1-(3/2)\delta$,
%
\begin{eqnarray}
\label{eqfirstbound1} L(\tilde{f}) \leq L^* + C \biggl(\frac{\log(\cN
_2(\F,\varepsilon,
S)/\delta)}{n}+ \max
_{i=1,\ldots,N}\hat{\Rad}_n\bigl(\ell\circ\hat{
\F}_i^S,S'\bigr) \biggr) .
\end{eqnarray}
%

%
To complete the proof of (\ref{2}) we need to evaluate the Rademacher
complexities appearing in~(\ref{eqfirstbound1}):
\[
\hat{\Rad}_n\bigl(\ell\circ\hat{\F}_i^S,S'
\bigr) = \En_{\sigma} \biggl[ \sup_{f\in\hat{\F}_i^S}
\frac{1}{n}\sum_{(x,y)\in S'} \sigma_i
\bigl(f(x)-y\bigr)^2 \biggr] . %
\]
The difficulty here is that the set $\hat{\F}_i^{S}=\hat{\F
}_i^{S}(\varepsilon)$
is defined via the pseudo-metric $d_S$ based on sample $S$ while the
empirical Rademacher complexity is evaluated on another sample $S'$. To
match the metrics, we embed $\hat{\F}_i^{S}(\varepsilon)$ into
$d_{S'}$-balls with properly chosen radius $\bar\gamma$:
%
\[
\hat{\F}^{S,S'}_i(\bar\gamma) \triangleq\bigl\{f\in\F\dvt
d_{S'}(f,\hat{c}_i) \leq\bar\gamma\bigr\}, \qquad i=1,\ldots,N,
\]
where the pseudo-metric $d_{S'}$ is taken with respect to the set
$S'$ while the $\varepsilon$-net $\hat{c}_1,\ldots,\hat{c}_N$ is
constructed with respect to $d_S$. The next lemma shows that, with
high probability, $\hat{\F}_i^{S}(\varepsilon)$ is included into
$\hat{\F}^{S,S'}_i(\bar\gamma)$ for an appropriate choice of
$\bar\gamma$.

%
\begin{lemma}
\label{leminclusion}
Fix $t>0$, $\varepsilon>0$. Let $r^*=r^*(\G)$ for $\G=\{(f-g)^2\dvt
f,g\in
\F\}$. Define $r_0 = (t+6\log\log n)/n$ and $\bar\gamma= \sqrt
{4\varepsilon^2 + 284r^* + 120 r_0}$. Then, with probability at least
$1-8N\mathrm{e}^{-t}$ with respect to the distribution of $S\cup S'$, we have
the inclusions
\[
\hat{\F}^S_i(\varepsilon) \subseteq\hat{
\F}^{S,S'}_i(\bar\gamma), \qquad i=1,\ldots,N,
\]
and hence, with the same probability,
\begin{eqnarray*}
&&\hat{\Rad}_n \bigl(\ell\circ\hat{\F}^{S}_i(
\varepsilon) , S' \bigr) \leq\hat{\Rad}_n \bigl(\ell\circ
\hat{\F}^{S,S'}_i(\bar\gamma), S' \bigr) ,
\qquad i=1,\ldots,N.
\end{eqnarray*}
%
\end{lemma}

\begin{pf} 
Let $P_n$ and $P'_n$ denote the empirical averages over the samples $S$
and $S'$, respectively. By Theorem~\ref{thmlocalization},
%
with probability at least $1-4\mathrm{e}^{-t}$,
\[
P(f-g)^2 \leq2P_n (f-g)^2 + 106 r^* + 48
r_0 \qquad\forall f,g\in\F, %
\]
and, with the same probability,
\[
P_n' (f-g)^2 \leq2P (f-g)^2 +
72 r^* + 24 r_0 \qquad\forall f,g\in\F.
\]
Therefore, with probability at least $1-8\mathrm{e}^{-t}$,
\[
P'_n(f-g)^2 \leq4 P_n(f-g)^2
+ 284 r^* + 120 r_0 \qquad\forall f,g\in\F.
\]
Applying this to $g=\hat{c}_i$ and taking a union bound over
$i=1,\ldots,N$, completes the proof.
\end{pf}
%

The next lemma gives an upper bound on the Rademacher complexity of
the set $\ell\circ\hat{\F}^{S,S'}_i(\bar\gamma)$.
%
\begin{lemma}
\label{lemradupperbound}
Let $r^*=r^*(\G)$ for $\G=\{(f-g)^2\dvt f,g\in\F\}$. Then, for any
$\bar\gamma\geq\sqrt{r^*}$ we have
\begin{eqnarray*}
\hat{\Rad}_n \bigl(\ell\circ\hat{\F}^{S,S'}_i(
\bar\gamma), S' \bigr) \leq\bar\gamma\sqrt{r^*} + \inf
_{\alpha\geq0} \biggl\{ 4\alpha+ \frac{24}{\sqrt{n}}\int
_{\alpha}^{\bar\gamma} \sqrt{\log\cN_2\bigl(\F,
\rho, S'\bigr)} \,\mathrm{d}\rho\biggr\} .
\end{eqnarray*}
\end{lemma}
\begin{pf}
Throughout the proof, we fix the samples $S$ and $S'$. We have
%
\begin{eqnarray}
\label{eqradsplittwo} \hat{\Rad}_n \bigl(\ell\circ\hat{
\F}^{S,S'}_i(\bar\gamma), S' \bigr)
&=& \En_\sigma\sup_{f\in\hat{\F}^{S,S'}_i(\bar\gamma)} \frac
{1}{n}
\sum_{(x_j,y_j)\in S'} \sigma_j
\bigl(f(x_j)-y_j\bigr)^2
\nonumber
\\
&=&\En_\sigma\sup
_{f\in\hat{\F}^{S,S'}_i(\bar\gamma)} \frac
{1}{n} \sum_{(x_j,y_j)\in S'}
\sigma_j \bigl(f(x_j)-\hat{c}_i(x_j)
\bigr)^2
\\
&&{}+ 2\En_\sigma\sup_{f\in\hat{\F}^{S,S'}_i(\bar
\gamma)} \frac{1}{n} \sum
_{(x_j,y_j)\in S'}\sigma_j\bigl(f(x_j)-
\hat{c}_i(x_j)\bigr) \bigl(\hat{c}_i(x_j)-y_j
\bigr),\nonumber
\end{eqnarray}
where we have used the decomposition $(f(x)-y)^2=(f(x)-\hat
{c}_i(x))^2+ (\hat{c}_i(x)-y)^2 + 2(f(x)-\hat{c}_i(x))(\hat
{c}_i(x)-y)$, $\forall x,y$, and the fact that $(\hat{c}_i(x)-y)^2$
does not depend on $f$. Conditionally on the sample $S$, the functions
$\hat{c}_i$ are fixed. Consider the sets of functions
\begin{eqnarray*}
\G'_i&=& \bigl\{(f-\hat{c}_i)^2\dvt
f\in\hat{\F}^{S,S'}_i(\bar\gamma) \bigr\}
\\
&=& \biggl\{(f-\hat{c}_i)^2\dvt f\in\F,
\frac{1}{n}\sum_{(x,y)\in
S'} \bigl(f(x)-
\hat{c}_i(x)\bigr)^2 \le\bar\gamma^2 \biggr\}
.
\end{eqnarray*}
Recall that we assume $\hat{c}_i\in\F$ (the $\varepsilon$-net is
proper). Thus $\G'_i\subseteq\G[\bar\gamma^2,S']$ for $\G=\{
(f-g)^2\dvt f,g\in\F\}$, which implies 
%
\begin{eqnarray}
\label{eqradsplittwo1} \En_\sigma\sup_{f\in\hat{\F}^{S,S'}_i(\bar\gamma)}
\frac{1}{n} \sum_{(x_j,y_j)\in S'} \sigma_j
\bigl(f(x_j)-\hat{c}_i(x_j)
\bigr)^2 &\leq&\hat{\Rad}_n\bigl(\G\bigl[\bar
\gamma^2, S'\bigr], S'\bigr)
\nonumber\\[-8pt]\\[-8pt]\nonumber
&\leq&\phi_n\bigl(\bar\gamma^2\bigr)\leq\bar\gamma
\sqrt{r^*},
\end{eqnarray}
where $\phi_n(\bar\gamma^2)=\phi_n(\bar\gamma^2,\G)$ and the last
inequality is due to the assumption $\bar\gamma^2> r^*$ and the
fact that $\phi_n(r)/\sqrt{r}$ is non-increasing.

We now turn to the cross-product term in (\ref{eqradsplittwo}).
Define the following sets of functions on $\X\times\Y$:
\[
\G_i^{S,S'} = \bigl\{g_f(x,y) = \bigl(f(x)-
\hat{c}_i(x)\bigr) \bigl(\hat{c}_i(x)-y\bigr)\dvt f\in
\hat{\F}^{S,S'}_i(\bar\gamma) \bigr\} . %
\]
Then,
%
\begin{eqnarray}
\label{eqradcrossterm1} \En_\sigma\sup_{f\in\hat{\F}^{S,S'}_i(\bar\gamma)}
\frac{1}{n} \sum_{(x_j,y_j)\in S'}\sigma_j
\bigl(f(x_j)-\hat{c}_i(x_j)\bigr) \bigl(
\hat{c}_i(x_j)-y_j\bigr) & =&\hat{
\Rad}_n\bigl(\G_i^{S,S'}, S'\bigr)
.
\end{eqnarray}
Observe that, for any $g_f\in\G_i^{S,S'}$,
%
\begin{eqnarray}
\label{eqradcrossterm0} \frac{1}{n}\sum_{(x,y)\in S'}
g_f(x,y)^2 &=& \frac{1}{n}\sum
_{(x,y)\in S'} \bigl(f(x)-\hat{c}_i(x)
\bigr)^2\bigl(\hat{c}_i(x)-y\bigr)^2
\\
&\leq&\frac{1}{n}\sum_{(x,y)\in S'} \bigl(f(x)-
\hat{c}_i(x)\bigr)^2\le\bar\gamma^2
\end{eqnarray}
since $\hat{c}_i$ and $y$ take values in $\Y=[0,1]$. For the same reason,
\begin{eqnarray*}
\frac{1}{n}\sum_{(x,y)\in S'} \bigl(g_f(x,y)-g_h(x,y)
\bigr)^2 &=&\frac
{1}{n}\sum_{(x,y)\in S'}
\bigl(f(x)-h(x)\bigr)^2\bigl(\hat{c}_i(x)-y
\bigr)^2
\\
&\leq&\frac{1}{n}\sum_{(x,y)\in S'} \bigl(f(x)-h(x)
\bigr)^2
\end{eqnarray*}
implying $\cN_2(\G^{S,S'}_i, \rho, S')\leq
\cN_2(\hat{\F}^{S,S'}_i(\bar\gamma), \rho, S')$ for all $\rho>0$.
Hence, by Lemma~\ref{lemradent},
%
\begin{eqnarray}
\label{eqradcrossterm} \hat{\Rad}_n\bigl(\G_i^{S,S'},
S'\bigr) &\leq&\inf_{\alpha\geq0} \biggl\{ 4\alpha+
\frac{12}{\sqrt{n}}\int_{\alpha}^{\bar\gamma} \sqrt{
\log\cN_2\bigl(\hat{\F}^{S,S'}_i(\bar\gamma),
\rho, S'\bigr)} \,\mathrm{d}\rho\biggr\}
\nonumber\\[-8pt]\\[-8pt]\nonumber
&\leq&\inf_{\alpha\geq0} \biggl\{ 4\alpha+ \frac{12}{\sqrt{n}}\int
_{\alpha}^{\bar\gamma} \sqrt{\log\cN_2\bigl(\F,
\rho, S'\bigr)} \,\mathrm{d}\rho\biggr\},
\end{eqnarray}
where the integration goes to $\bar\gamma$ in view of
(\ref{eqradcrossterm0}). The lemma now follows from
(\ref{eqradsplittwo})--(\ref{eqradcrossterm1}) and~(\ref{eqradcrossterm}).
\end{pf}

Combining (\ref{eqfirstbound1}), Lemma~\ref{leminclusion} with
$t=\log(16 N/\delta)$, and Lemma~\ref{lemradupperbound} we find
that, with probability at least $1-2\delta$,
%
\begin{eqnarray}
\label{eqboundfin} L(\tilde{f}) &\leq& L^* + C \biggl(\frac{\log(\cN_2(\F
,\varepsilon, S)/\delta)}{n}+ \bar\gamma
\sqrt{r^*}
\nonumber\\[-8pt]\\[-8pt]\nonumber
&&\hspace*{36pt}{} + \inf_{\alpha\geq0} \biggl\{ \alpha+ \frac{1}{\sqrt{n}}\int
_{\alpha}^{\bar\gamma} \sqrt{\log\cN_2\bigl(\F,
\rho, S'\bigr)} \,\mathrm{d}\rho\biggr\} \biggr),
\end{eqnarray}
which yields the bound (\ref{1}).

\section{Proofs of the lower bounds}\label{secproofslower}

\begin{pf*}{Proof of Theorem~\ref{thlow1}}
Fix some $0<\alpha<1$ and set $k=\lceil d/\alpha\rceil$.
Let $\cC$ be the set of all binary sequences $\omega\in\{0,1\}^k$
with at most $d$ non-zero components. By the $d$-selection lemma (see,
e.g., Lemma 4 in \cite{RagRak11b}), for $k\geq2d$ there exists of a
subset $\cC'$ of $\cC$ with the following properties: (a) $\log\llvert
\cC
'\rrvert\geq(d/4)\log(k/(6d))$ and (b) $\rho_H(\omega, \omega
')\ge d$
for any $\omega, \omega'\in\cC'$. Here, $\rho_H(\omega, \omega
')=\sum_{j}{\mathbf1} \{\omega_j \ne\omega'_j \}$ denotes
the Hamming
distance where $\omega_j,\omega'_j$ are the components of $\omega
,\omega'$. To any $\omega\in\cC'$ we associate a function $g_\omega
$ on $\X$ defined by $g_\omega(x^i)=\omega_i$ for $ i=1,\ldots, k$
and $g_\omega(x^i)=0$, $i \ge k+1$, where $\omega_i$ is the $i$th
component of $\omega$.

Let $\mu_X$ be the distribution on $\X$ which is uniform on $\{
x^1,\ldots,x^k\}$, putting probability $1/k$ on each of these $x^j$
and probability 0 on all $x^j$ with $j\ge k+1$. Denote by $\P_\omega$
the joint distribution of $(X,Y)$ having this marginal $\mu_X$ and
$Y\in\{0,1\}$ with the conditional distribution $\En(Y\llvert X=x)=
P(Y=1\rrvert X=x) = 1/2+g_\omega(x)/4 \triangleq\eta_\omega
(x)$ for all
$x\in\X$.

Consider now a set of functions
$\F'=\{\eta_\omega\dvt \omega\in\cC' \} \subset\F$.
Observe that, by construction,
%
\begin{equation}
\label{lb1}\llVert\eta_{\omega}-\eta_{\omega
'}\rrVert
^2 = \rho_H\bigl(\omega,\omega'\bigr)
/(16k) \geq\alpha/32 \qquad\forall\omega,\omega'\in
\cC'.
\end{equation}
On the other hand, the Kullback--Leibler divergence between $\P_\omega
$ and $\P_{\omega'}$ has the form
\[
K(\P_\omega,\P_{\omega'})=n \En\biggl(\eta_\omega(X)
\log\frac
{\eta_\omega(X)}{\eta_{\omega'}(X)}+\bigl(1-\eta_\omega(X)\bigr)\log
\frac
{(1-\eta_\omega(X))}{(1-\eta_{\omega'}(X))}
\biggr) . %
\]
Using the inequality $-\log(1+u)\le-u+u^2/2$, $\forall u>-1$, and the
fact that $1/2\le\eta_\omega(X) \le3/4$ for all $\omega\in\cC'$
we obtain that the expression under the expectation in the previous
display is bounded by $2( \eta_\omega(X)-\eta_{\omega'}(X))^2$,
which implies
%
\begin{equation}
\label{lb2} K(\P_\omega,\P_{\omega'})\le2n \En\bigl(
\eta_\omega(X)-\eta_{\omega
'}(X)\bigr)^2 \le
\frac{n\llVert g_\omega-g_{\omega'}\rrVert
^2}{8}\le\frac
{nd}{8k}\le\frac{n\alpha}{8} \qquad\forall\omega,
\omega'\in\cC'.
\end{equation}
From (\ref{lb1}), (\ref{lb2}) and Theorem 2.7 in \cite{Tsy09}, the
result of Theorem \ref{thlow1} follows if we show that
%
\begin{equation}
\label{lb3} n\alpha/{8} \le\log\bigl(\bigl\llvert\F'\bigr
\rrvert-1\bigr)/16
\end{equation}
with
\[
\alpha= C_1\frac{d}{n}\log\frac{C_2n}{d},
\]
where $C_1,C_2>0$ are constants. Assume first that $d\ge4$. Then,
using the inequalities $ \log(\llvert\F'\rrvert -1)\ge
\log(\llvert\cC'\rrvert/2 \geq
(d/4)\log(k/(6d)) -\log2 \ge(d/4)\log(1/(12\alpha))$
it is enough to show that
\[
n\alpha\le\frac{d}{8}\log\frac{1}{12\alpha} . %
\]
Using that $x\ge2\log x$ for $x\ge0$ it is easy to check that the
inequality in the last display holds if we choose, for example,
$C_1=1/16, C_2=1/(12C_1)$.
In the case $d\le3$, it is enough to consider $\alpha=(C_1/n)\log
(C_2n)$ and (\ref{lb3}) is also satisfied for suitable $C_1,C_2$.
\end{pf*}

\begin{pf*}{Proof of Theorem~\ref{lemlowerrich}}
We first prove the entropy bound (\ref{thlow21}). It suffices to
obtain the same bound for $B_p$ in place of $\F$. Fix $\varepsilon>0$
and set $J=(2/\varepsilon)^p$. Without loss of generality, assume
that $J$ is an integer. Let $\M$ be an $\varepsilon$-net on $B_p$ in
$\ell_\infty$ metric constructed as follows. For all $v\in\M$ , the
coordinate $v_j$ of $v$ takes discrete values with step $\varepsilon$
within the interval $[-j^{-1/p},j^{-1/p}]$ if $j\le J$, and $v_j=0$ for
all $j>J$. Then,
\[
\llvert\M\rrvert\le\prod_{j=1}^J
\biggl(\frac
{2}{\varepsilon j^{1/p}} \biggr).
\]
One can check that
\[
\log\Biggl(\prod_{j=1}^J
j^{-1/p} \Biggr)= -\frac{1}{p}\sum_{j=1}^J
\log j\le-\frac{1}{p}\int_2^J (\log t)
\,\mathrm{d}t \le-\frac
{J}{p}(\log J-1),
\]
which implies that $\llvert\M\rrvert \le\exp(J/p)$. Thus
(\ref{thlow21}) follows.

\textit{Proof of} (\ref{thlow22}).
%
Fix $d=\lceil n^{p/(2+p)}\rceil$. Let $\Omega_d=\{0,1\}^d$ be the
set of all binary sequences of length~$d$. Define $\mu_X$ as the
distribution on $\X$ which is uniform on $\{\e_1,\ldots,\e_d\}$,
putting probability $1/d$ on each of these $\e_j$ and probability 0
on all $\e_j$ with $j\ge d+1$. For any $\omega\in\Omega_d$, denote
by $\P_\omega$ the joint distribution of $(X,Y)$ having this
marginal $\mu_X$ and $Y\in\{0,1\}$ with the conditional distribution
defined by the relation
\[
\eta_\omega(\e_i)\triangleq\En(Y \mid X =
\e_i) = P(Y=1\mid X=\e_i) = \frac{1}{2} +
\frac{\omega_i}{4d^{1/p}}
\]
for\vspace*{1pt} $i=1,\ldots,d$, and $\eta_\omega(\e_i)=1/2$ for $i\ge d+1$. The
regression function corresponding to $\P_\omega$ is then
$\eta_\omega= \{\eta_\omega(\e_j)\} \in\ell$.
It is easy to see that since $\omega_i\in\{0,1\}$ for any
estimator $\hat{f}=\{\hat{f}_j\}\in\ell$ we have
\[
\bigl\llvert\hat{f}_i-\eta_\omega(\e_i)\bigr
\rrvert\ge\frac
{1}{2}\biggl\llvert\frac{1}{2} +\frac{\hat{\omega}_i}{4d^{1/p}}-
\eta_\omega(\e_i)\biggr\rrvert= \frac{\llvert\hat{\omega}_i-\omega
_i\rrvert}{8d^{1/p}}, \qquad
i=1,\ldots,d,
\]
where $\hat{\omega}_i$ is the closest to $4d^{1/p}(\hat{f}_i-1/2)$
element of the set $\{0, 1\}$. Therefore,
%
\begin{eqnarray}
\label{thlow24} \llVert\hat{f}-\eta_\omega\rrVert^2 &\ge&
\frac{1}{d} \sum_{i=1}^d
\frac
{\llvert\hat{\omega}_i-\omega_i\rrvert^2}{64
d^{2/p}}=\frac{\rho_H(\hat
{\omega},\omega)}{64 d^{1+2/p}},
\end{eqnarray}
where $\rho_H(\cdot,\cdot)$ is the Hamming distance. From Assouad's
lemma (cf. Theorem 2.12(iv) in \cite{Tsy09}),
%
\begin{eqnarray}
\label{thlow25} \inf_{\hat{\omega}}\max_{\omega\in\Omega_d}\mathbf
{E}^{(n)}_\omega\rho_H(\hat{\omega},\omega) \ge
\frac{d}{4}\exp(-\alpha),
\end{eqnarray}
where $\alpha= \max\{K( \P_\omega, \P_{\omega'})\dvt
\omega,\omega'\in\Omega_d, \rho_H(\omega,\omega')=1\}$. Here,
$\mathbf{E}^{(n)}_\omega$ denotes
the distribution of the $n$-sample $D_n$ when $(X_i,Y_i)\sim
\P_\omega$ for all $i$. Since $1/2\le\eta_\omega(X)\le3/4$, the
Kullback--Leibler divergence can be bounded in the same way as in (\ref{lb2}):
\begin{eqnarray*}
K( \P_\omega, \P_{\omega'}) &\le&2n \En\bigl(\eta_\omega(X)-
\eta_{\omega'}(X)\bigr)^2 = \frac{2n}{d} \sum
_{i=1}^d \frac{(\omega
_i-\omega'_i)^2}{64 d^{2/p}}=\frac{n\rho_H(\hat{\omega},\omega
)}{32 d^{1+2/p}} \le
\frac{1}{32}
\end{eqnarray*}
for all $\omega,\omega'\in\Omega_d$ such that
$\rho_H(\omega,\omega')=1$. Combining this result with
(\ref{thlow24}) and (\ref{thlow25}), we find
%
\begin{eqnarray}
\label{thlow26} \inf_{\hat{f}}\max_{\omega\in\Omega_d}
\mathbf{E}^{(n)}_\omega\llVert\hat{f}-\eta_\omega
\rrVert^2 \ge\frac{\mathrm{e}^{-1/32}}{128 d^{2/p}} \ge c_*n^{-2/(2+p)}
\end{eqnarray}
for some absolute constant $c_*>0$. Now, the set $\{\eta_\omega\dvt
\omega\in\Omega_d\}$ is contained in $\F$, so that
%
\begin{eqnarray}
\label{thlow26a} W_n(\F)\ge\inf_{\hat{f}}\max
_{\omega\in\Omega_d} \mathbf{E}^{(n)}_\omega\llVert
\hat{f}-\eta_\omega\rrVert^2
\end{eqnarray}
and (\ref{thlow22}) follows
immediately from (\ref{thlow26}) and (\ref{thlow26a}).

\textit{Proof of} (\ref{thlow23}).
Set $d=2\lceil n^{p/(p-1)}\rceil$
and define the
joint distribution $\P_\omega$ of $(X,Y)$ as in the proof of
(\ref{thlow22}) with the difference that now we
choose the conditional probabilities as follows:
\[
\eta_\omega(\e_j)=\frac{1}{2}+\frac{\omega_j}{4} ,
\qquad j=1,\ldots,d\quad\mbox{and}\quad\eta_\omega(\e_j) =
\frac{1}{2}, \qquad j\ge d+1,
\]
where $\omega=(\omega_1,\ldots,\omega_d)\in\Omega_d'=\{-1,1\}^d$.
Set $\eta_\omega=
\{\eta_\omega(\e_j)\}\in\ell$ with $\omega\in\Omega_d'$. Then
\[
\inf_{f\in\F} \llVert f-\eta_\omega\rrVert
^2 \le\llVert f_\omega-\eta_\omega\rrVert
^2= \frac{1}{16}\bigl(1-d^{-1/p}\bigr)^2,
\]
where $f_\omega= \{f_\omega(\e_j)\}\in\F$ is
a sequence with components
\[
f_\omega(\e_j) = \frac{1}{2}+\frac{\omega_j}{4d^{1/p}} ,
\qquad j=1,\ldots,n\quad\mbox{and}\quad f_\omega(\e_j) =
\frac{1}{2}, \qquad j\ge d+1.
\]
Hence,
%
\begin{eqnarray*}
V_n(\F) &=& \inf_{\hat{f}}\sup_{P_{XY}\in\cP}
\Bigl\{ \En\llVert\hat{f}-\eta\rrVert^2 -\inf_{f\in\F}
\llVert f-\eta\rrVert^2 \Bigr\}
\\
&\ge&\inf_{\hat{f}}\max_{\omega\in\Omega_d'} \Bigl\{
\mathbf{E}^{(n)} _\omega\llVert\hat{f}-\eta_\omega
\rrVert^2 -\inf_{f\in
\F} \llVert f-
\eta_\omega\rrVert^2 \Bigr\}
\\
& \ge& \inf_{\hat{f}} \int_{\Omega_d'}
\mathbf{E}^{(n)}_\omega\llVert\hat{f}-\eta_\omega
\rrVert^2 \nu(\mathrm{d}\omega) - \frac{1}{16}\bigl(1-d^{-1/p}
\bigr)^2,
\end{eqnarray*}
where $\nu$ is the probability measure on $\Omega_d'$ under which
$\omega_1,\ldots,\omega_d$ are i.i.d. Rademacher random variables.
Passing to sequences $\hspace*{4pt}{\bar{\hspace*{-4pt}\hat f}}$, ${\bar\eta}_\omega$ in $\ell
$ with components $\hspace*{4pt}{\bar{\hspace*{-4pt}\hat f}}_j={\hat
f}_j-1/2$, ${\bar\eta}_\omega(\e_j)=\eta_\omega(\e_j)-1/2$,
respectively, we may
write
\begin{eqnarray*}
V_n(\F) &\ge&\inf_{\hspace*{4pt}{\bar{\hspace*{-4pt}\hat f}}} \int_{\Omega_d'}
\mathbf{E}^{(n)}_\omega\llVert\hspace*{4pt}{\bar{\hspace*{-4pt}\hat f}}-{\bar
\eta}_\omega\rrVert^2 \nu(\mathrm{d}\omega) - \frac{1}{16}
\bigl(1-d^{-1/p}\bigr)^2.
\end{eqnarray*}
For $j=1,\ldots,d$, denote by $\hat f[j]$ and $r_j$ the components of
$\hspace*{4pt}{\bar{\hspace*{-4pt}\hat f}}$ and of ${\bar\eta}_\omega$, respectively. We will
sometimes write $\hat f[j]=\hat f[j,D_n]$ to emphasize the dependence
on the sample $D_n = \{(X_1,Y_1),\ldots,(X_n,Y_n)\}$. Then, we can rewrite
the above integral in the form
\[
\int_{\Omega_d'} \mathbf{E}^{(n)}_\omega\llVert
\hspace*{4pt}{\bar{\hspace*{-4pt}\hat f}}-{\bar\eta}_\omega\rrVert^2 \nu(\mathrm{d}\omega)= {
\En}_{r_1,\ldots, r_d} {\En}_{D_n} \Biggl[\frac{1}d \sum
_{j=1}^d\bigl(r_j-\hat f[j]
\bigr)^2 \Biggr],
\]
where ${\En}_{r_1,\ldots, r_d}$ and ${\En}_{D_n}$ denote the
expectation over the joint distribution of $r_1,\ldots, r_d$ and
over the distribution of $D_n$ given $r_1,\ldots, r_d$, respectively.

Consider the random vector composed of indicators $\zeta=(I(X_1=e_j),
\dots, I(X_n=e_j))$. For any $j$ and any fixed $r_1,\ldots,r_d$,
\begin{eqnarray*}
{\En}_{D_n} \bigl[\bigl(r_j-\hat f[j,D_n]
\bigr)^2\bigr]&=&{\En}_{\zeta} {\En}_{D_n} \bigl[
\bigl(r_j-\hat f[j,D_n]\bigr)^2\mid\zeta
\bigr]
\\
&\ge& {\mathbf P}(\zeta=0) {\En}_{D_n} \bigl[\bigl(r_j-
\hat f[j,D_n]\bigr)^2\mid\zeta=0\bigr]
\\
&\ge& {\mathbf P}(\zeta=0) \bigl(r_j-{\En}_{D_{n}}\bigl[
\hat f[j,D_n]\mid\zeta=0\bigr] \bigr)^2,
\end{eqnarray*}
where we have used Jensen's inequality. We may write ${\En
}_{D_{n}}[\hat f[j,D_n]\mid\zeta=0]$ in the form
\[
{\En}_{D_{n}}\bigl[\hat f[j,D_n]\mid\zeta=0\bigr]=G\bigl(
\{r_k\dvt k\ne j\}\bigr),
\]
where $G$ is some measurable function. Indeed, under the condition
$\zeta=0$ the distribution of $D_{n}$ coincides with that of $\{
(X_i,Y_i)\dvt X_i\ne e_j\}$, which is entirely defined by $\{r_k\dvt
k\ne
j\}$. Thus,
\begin{eqnarray*}
&&{\En}_{r_j} {\En}_{D_n} \bigl[\bigl(r_j-\hat
f[j,D_n]\bigr)^2\bigr]
\\
&&\quad \ge{\mathbf P} (\zeta=0) {
\En}_{r_j}\bigl[\bigl(r_j-G\bigl(\{r_k\dvt k\ne j
\}\bigr)\bigr)^2\bigr]
\\
&&\quad = {\mathbf P}(\zeta=0) \biggl[\frac{1}2 \biggl(\frac{1}4-G
\bigl(\{r_k\dvt k\ne j\}\bigr) \biggr)^2 +
\frac{1}2 \biggl(-\frac{1}4-G\bigl(\{r_k\dvt k\ne j\}
\bigr) \biggr)^2 \biggr]
\\
&&\quad \ge\frac{1}{16}{\mathbf P}(\zeta=0)= \frac{1}{16} \biggl(1-
\frac
{1}{d} \biggr)^n,
\end{eqnarray*}
where ${\En}_{r_j}$ denotes the expectation over the distribution of
$r_j$ and we have used that $r_j$ takes values $1/4$ and $-1/4$ with
probabilities $1/2$. This implies
\[
\inf_{\hat f} {\En}_{r_1,\ldots, r_d} {\En}_{D_n}
\Biggl[\frac{1}d \sum_{j=1}^d
\bigl(r_j-\hat f[j]\bigr)^2 \Biggr]\ge\frac{1}{16}
\biggl(1-\frac
{1}{d} \biggr)^n, %
\]
so that
\[
V_n(\F) \ge\frac{1}{16} \biggl[ \biggl(1-\frac{1}{d}
\biggr)^n - \bigl(1-d^{-1/p}\bigr)^2 \biggr].
\]
Using that $1-x\ge\exp(-3x/2)$ for $0<x\le1/2$ we have $
(1-\frac{1}{d} )^n\ge\exp(-3n/(2d))\ge1- 3n/(2d)$ for $d\ge
2n$. Since $d =2\lceil n^{p/(p-1)}\rceil$ we find
\begin{eqnarray*}
V_n(\F) &\ge& 1- 3n/(2d) - \bigl(1-d^{-1/p}
\bigr)^2= - 3n/(2d) + 2d^{-1/p} - d^{-2/p}
\\
&\ge&- 3n/(2d) - d^{-1/p}\ge d^{-1/p}/4\ge cn^{-1/(p-1)}
\end{eqnarray*}
for some absolute constant $c>0$.
\end{pf*}

\setcounter{equation}{0}
\begin{appendix}\label{secappendix}
\section*{Appendix}

The following result is a modification of Theorem 6.1 in \cite
{bousquet2002concentration}. 

%
\begin{theorem}
\label{thmlocalization}
Let $\G$ be a class of non-negative functions bounded by $b$ and
admitting a localization radius $r^*=r^*(\G)$. 
Then for all $n\ge5$ and $t>0$, with probability at least $1-4\mathrm{e}^{-t}$,
for all $g\in\G$ we have
%
\begin{eqnarray}
\label{eqthmloc1}  Pg &\leq&2 P_n g + 106 r^* + 48 r_0,
\\
\label{eqthmloc2}  P_n g &\leq&2 P g
+ 72r^* + 24 r_0,
\end{eqnarray}
%
where $r_0 = b(t+6\log\log n)/n$.
\end{theorem}
\begin{pf*}{Proof of Theorem~\ref{thmlocalization}}
The fact that for $n\ge5$ inequality~(\ref{eqthmloc1}) holds with
probability at least $1-\mathrm{e}^{-t}$ for all $g\in\G$ is proved in Theorem
6.1 in \cite{bousquet2002concentration}. Moreover, it is shown in the
proof of that theorem that, on the same event of probability at least
$1-\mathrm{e}^{-t}$ (denote this event by ${\mathcal B}$),
%
\begin{eqnarray}
\label{eqbou1} Pg \leq P_n g + \sqrt{Pg}\bigl(\sqrt{8r^*} +
\sqrt{4r_0}\bigr) + 45 r^* + 20 r_0\qquad\forall g\in
\bigcup_{k=0}^{k_0}\G_k,
\end{eqnarray}
where $\G_k = \{g\in\G\dvt \delta_{k+1}\leq Pg \leq\delta_k\}$,
$\delta_k = b2^{-k}$ for $k\geq0$, and $k_0>0$ be the largest integer
such that $\delta_{k_0+1} \geq b/n$. A straightforward modification of
the argument in~\cite{bousquet2002concentration} leading to~(\ref
{eqbou1}) yields that, on the event ${\mathcal B}$,
%
\begin{eqnarray}
\label{eqbou1a} \llvert Pg - P_n g\rrvert\le\sqrt{Pg}\bigl(
\sqrt{8r^*} + \sqrt{4r_0}\bigr) + 45 r^* + 20 r_0\qquad
\forall g\in\bigcup_{k=0}^{k_0}
\G_k,
\end{eqnarray}
so that
%
\begin{eqnarray}\label{eqbou2}
P_n g& \leq& Pg + \sqrt{Pg}\bigl(\sqrt{8r^*} +
\sqrt{4r_0}\bigr) + 45 r^* + 20 r_0
\nonumber\\[-8pt]\\[-8pt]\nonumber
& \leq&2Pg + 53 r^* + 24 r_0\qquad\forall g\in\bigcup
_{k=0}^{k_0}\G_k,
\end{eqnarray}
proving~(\ref{eqthmloc2}) for $g\in\bigcup_{k=0}^{k_0}\G_k$ with
probability at least $1-\mathrm{e}^{-t}$.

Now, consider $g\in\G_*=\G\setminus\bigcup_{k=0}^{k_0}\G_k$. First,
for any $g\in\G_*$, $P g \leq\delta_k \leq\delta_{k_0}\leq4b/n$.
Hence $\G_*\subseteq\G'=\{g\in\G\dvt Pg < 4b/n\}$. By Lemma 6.1 in
\cite{bousquet2002concentration}, with probability at least $1-3\mathrm{e}^{-t}$,
%
\begin{eqnarray}
\label{eqbou3} \llvert P_n g-P g\rrvert&\leq& 6 {\hat
\Rad}_n\bigl(\G',S\bigr) + \frac
{b}{n}(
\sqrt{2t}+6t) \nonumber\\[-8pt]\\[-8pt]\nonumber
&\leq& 6 {\hat\Rad}_n\bigl(\G',S\bigr) +
\frac{b(7t+1)}{n} \qquad\forall g\in\G'.\nonumber
\end{eqnarray}
Denote the event where~(\ref{eqbou3}) holds by ${\mathcal B}'$, and
define
\[
U' = 6 {\hat\Rad}_n\bigl(\G',S\bigr) +
Pg + \frac{b(7t+1)}{n} .
\]
On the event ${\mathcal B}'$ we have $P_n g \leq U'$ for any
$g\in\G'$, so that
\[
{\hat\Rad}_n\bigl(\G',S\bigr)\leq{\hat
\Rad}_n\bigl(\bigl\{g\in\G\dvt P_n g\leq U'
\bigr\},S\bigr) \leq\phi_n\bigl(U'\bigr),
\]
where $\phi_n(\cdot)=\phi_n(\cdot,\G)$ is an upper function for
$\G$ satisfying the sub-root property. In view of this property,
\[
U' \leq6\phi_n\bigl(U'\bigr) + Pg +
\frac{b(7t+1)}{n} \leq6\sqrt{U'}\sqrt{r^*} + Pg +
\frac{b(7t+1)}{n} .
\]
Solving for $\sqrt{U'}$ we get
\[
\sqrt{U'}\leq6\sqrt{r^*} + \sqrt{Pg + \frac{b(7t+1)}{n}}
\]
and thus, on the event ${\mathcal B}'$,
%
\begin{eqnarray}
\label{eqbou4} P_n g \leq U' \leq2Pg+ 72r^* +
\frac{2b(7t+1)}{n} \le2Pg+ 72r^* + 14r_0\qquad\forall g\in
\G',
\end{eqnarray}
where the last inequality is due to the fact that $7t+1\le7(t+6\log
\log n)$ for all $n\ge3$. Combining~(\ref{eqbou2}) and~(\ref
{eqbou4}), we then obtain~(\ref{eqthmloc2}) holds for all $g\in\G$
on the event ${\mathcal B}\cap{\mathcal B}'$ of probability at least
$1-4\mathrm{e}^{-t}$.
\end{pf*}

%
\begin{pf*}{Proof of Lemma~\ref{lemboundonr-star}}
\textit{Proof of} (i). We apply Lemma~2.2 in \cite{srebro2010smoothness}
for the loss function defined by $\varphi(t,y)=t^2, \forall t,y\in
\R$. The second derivative of this function with respect to the first
argument is $H=2$, that is, the function is 2-smooth in the terminology
of \cite{srebro2010smoothness}. Consider the class of differences $\cH
=\{f-g\dvt f,g\in\F\}$. Then Lemma~2.2 in \cite{srebro2010smoothness}
provides the following bound for the Rademacher complexity of the set
$\cL= \{ (x,y) \mapsto\varphi(h(x),y) \dvt h\in\cH, n^{-1} \sum
_{(x,y)\in S}h^2(x)\le r\}$:
\[
\hat{\Rad}_n(\cL,S)\leq21\sqrt{12r}\log^{3/2}(64n)
\Rad_n(\cH). %
\]
On the other hand, $\cL= \G[r,S]$, and $\Rad_n(\cH)\le2 \Rad_n(\F
)$, so that
%
\begin{eqnarray*}
\hat{\Rad}_n\bigl(\G[r,S],S\bigr) \leq42\sqrt{12r}
\log^{3/2}(64n) \Rad_n(\F).
\end{eqnarray*}
Now define the function $\phi_n(r)$ as the right-hand side of this
inequality. This immediately yields a localization radius
\[
r^* = 12\cdot42^2 \log^3(64 n) \Rad_n^2(
\F),
\]
and (\ref{eqrnstar-large}) follows.

\textit{Proof of} (ii). Let $(f-g)^2$ and $(\bar f-\bar g)^2$ be two
elements of $\G$, where $f,g,\bar f,\bar g\in\F$. Since all these
functions take values in $[0,1]$ we get that, for any $x\in\X$,
\[
\bigl(\bigl(f(x)-g(x)\bigr)^2- \bigl(\bar f(x)-\bar g(x)
\bigr)^2 \bigr)^2\le8 \bigl(\bigl(f(x)-\bar f(x)
\bigr)^2+\bigl(g(x)-\bar g(x)\bigr)^2 \bigr) . %
\]
Thus, if $d_S(f,\bar f)\le\varepsilon$ and $d_S(g,\bar g)\le\varepsilon$
for some $\varepsilon>0$, then $d_S((f-g)^2,(\bar f-\bar g)^2)\le
4\varepsilon$.
This implies the relation between the empirical entropies: $\cN_2(\G,
\rho, S) \leq\cN_2(\F,\rho/4,S)$ for all $\rho>0$. Using it
together with the bound $ \cN_2(\F,\rho/4,S)\le\max\{1, (4A/\rho
)^v\}$ and applying Lemma~\ref{lemradent} we obtain
\begin{eqnarray}
\label{equpperbdrademachercover} \hat{\Rad}_n\bigl(\G[r,S], S\bigr)
&\leq&
\frac{12}{\sqrt{n}}\int_{0}^{\sqrt
{r}} \sqrt{\log
\cN_2(\G, \rho, S)} \,\mathrm{d}\rho\nonumber
\\
&\leq& \frac{12}{\sqrt{n}}\int
_{0}^{\sqrt{r}\wedge1/(4A)} \sqrt{v\log(4A/\rho)} \,\mathrm{d}\rho
\\
&\leq&\frac{48 A\sqrt{v}}{\sqrt{n}}\int_{0}^{\sqrt{r}/(4A)\wedge
1} \sqrt{\log(1/t)} \,\mathrm{d}t\nonumber
\\
&\le& 24\sqrt{\frac{vr}{n}} \bigl(\log(4eA/\sqrt
{r})\vee1
\bigr)^{1/2},
\nonumber
\end{eqnarray}
where we have used that, integrating by parts,
%
\begin{eqnarray}
\label{eqvcintegral} \int_0^b \sqrt{\log(e/t)} \,\mathrm{d}t
&=&  b \sqrt{\log(e/b)} + (b/2) \bigl(\log(e/b)\bigr)^{-1/2}
\nonumber\\[-8pt]\\[-8pt]\nonumber
&\le& 2b \sqrt{
\log(e/b)}\qquad \forall0<b\le1.
\end{eqnarray}
In view of (\ref{equpperbdrademachercover}), we can take $\phi
_n(r)= 24\sqrt{\frac{vr}{n}} (\log(4eA/\sqrt{r})\vee1
)^{1/2}$ as an upper function in (\ref{eqphindef}).
Now, we are looking for $r^*$, which is an upper bound on the solution
of the equation $\phi_n(r)=r$. Since the function $u\mapsto(a/u)(\log
(b/u)\vee1)^{1/2} $, for $a, b>0$, is decreasing when $u>0$
one can check that $u^* = a(\log(b/a)\vee1)^{1/2}$ as an upper bound
on the solution of $(a/u)(\log(b/u)\vee1)^{1/2}=1$ whenever $b\ge e a>0$.
That is, for $n\ge Cv$ with $C>0$ large enough depending only on $A$,
we can take
%
\begin{eqnarray}
\label{eqrnvcclass} r^* = \biggl[24\sqrt{\frac{v}{n}} \biggl(\log\biggl(
\frac
{eA}{6}\sqrt{\frac{n}{v}} \biggr)\vee1 \biggr)^{1/2}
\biggr]^2 \le C\frac{v}{n}\log\biggl(\frac{en}{v} \biggr)
\end{eqnarray}
for some constant $C>0$ depending only on $A$.

\textit{Proof of} (iii). For a finite class $\F$, the covering numbers
satisfy $\cN_2(\F,\varepsilon, S) \leq\llvert\F\rrvert$
for all $\varepsilon>0$
and, along the lines of (\ref{equpperbdrademachercover}),
\[
\hat{\Rad}_n\bigl(\G[r,S], S\bigr) \leq\frac{12}{\sqrt{n}}\int
_{0}^{\sqrt
{r}} \sqrt{\log\cN_2(\F, \rho/4,
S)}\,\mathrm{d}\rho\le\frac{12\sqrt
{r\log\llvert\F\rrvert}}{\sqrt{n}}\triangleq\phi_n(r), %
\]
so that we can take
$r^* = 144 (\log\llvert\F\rrvert)/n$.
\end{pf*}

\begin{pf*}{Proof of Lemma~\ref{lemerm}}
Assume that there exists $f^*\in\F$ such that $L(f^*)=\min_{f\in\F
} L(f)$ (if this is not the case, an easy modification of the proof is
possible by considering an approximate minimizer). We apply Theorem~3.3
in \cite{bbm-lrc-05} to the class of functions $\G=\ell\circ\F
-\ell\circ f^*$. Observe that, for any $f\in\F$, the variance of the
random variable $\ell\circ f (X,Y)- \ell\circ f^*(X,Y)$ satisfies
\[
\operatorname{Var}\bigl(\ell\circ f - \ell\circ f^*\bigr) \leq\En\bigl[ \bigl(
\bigl(f(X)-Y\bigr)^2 - \bigl(f^*(X)-Y\bigr)^2
\bigr)^2\bigr]\leq2\bigl(L(f)-L\bigl(f^*\bigr)\bigr)
\]
and thus the assumption of Theorem 3.3 in \cite{bbm-lrc-05} holds with
$B=2$. Applying that theorem with $K=2$ we get that, for any $t>0$,
with probability at least $1-\mathrm{e}^{-t}$, for any $g\in\G$,
\begin{eqnarray*}
P g \leq2P_n g + c''_1 {\bar
r}^* + \frac{t(22 + c''_2)}{n},
\end{eqnarray*}
where $c''_1=704$, $c''_2=104$, and ${\bar r}^*$ is the solution of
fixed point equation $\psi(r)=r$, for a function $\psi$ satisfying
the sub-root property and the inequality $\psi(r) \ge2\En{\hat\Rad
}_n(\G\cap\{2Pg\le r\},S')$.
Choose now a constant function $\psi(r) \equiv2\En{\hat\Rad}_n(\G
,S')$, which trivially satisfies the sub-root property and has the
fixed point ${\bar r}^*=2\En{\hat\Rad}_n(\G,S')$. Since $ \En{\hat
\Rad}_n(\G,S')= \En{\hat\Rad}_n(\ell\circ\F,S')$, and $P_ng\le
0$ for $g= \ell\circ{\hat{f}^{\mathrm{emp}}}- \ell\circ{f^*}$, we
obtain that, with probability at least $1-\mathrm{e}^{-t}$,
\begin{eqnarray*}
L\bigl(\hat{f}^{\mathrm{emp}}\bigr) -L\bigl(f^*\bigr) \leq2c''_1
\En{\hat\Rad}_n\bigl(\ell\circ\F,S'\bigr) +
\frac{t(22 + c''_2)}{n} ,
\end{eqnarray*}
where we have used that $L(f)=P(\ell\circ f)$. Next, by Lemma A.4 in
\cite{bbm-lrc-05}, with probability at least $1-\mathrm{e}^{-t}$,
\[
\En{\hat\Rad}_n\bigl(\ell\circ\F,S'\bigr) \leq2{\hat
\Rad}_n\bigl(\ell\circ\F,S'\bigr) + \frac{t}{n}
.
\]
Combining the results of the last two displays we find that, with
probability at least $1-2\mathrm{e}^{-t}$,
\begin{eqnarray*}
L\bigl(\hat{f}^{\mathrm{emp}}\bigr) -L\bigl(f^*\bigr) \leq4c''_1
{\hat\Rad}_n\bigl(\ell\circ\F,S'\bigr) +
\frac{t(22 + c''_2+2c''_1)}{n} .
\end{eqnarray*}\upqed
\end{pf*}
\end{appendix}

\section*{Acknowledgements}
We gratefully acknowledge the support of NSF under Grants CAREER
DMS-0954737 and CCF-1116928. The work of the third author was supported
by GENES, and by the French National Research Agency (ANR) under the
Grants Idex ANR-11-IDEX-0003-02, Labex ECODEC (ANR-11-LABEX-0047), and
IPANEMA (ANR-13-BSH1-0004-02).



%

\printhistory
\end{document}